\documentclass[journal]{elsarticle}

\usepackage{hyperref}
\usepackage{epsfig,extarrows}
\usepackage{amsmath,amsthm}
\usepackage{amssymb}
\usepackage{graphicx,epstopdf,latexsym,amsfonts,fancyhdr}
\usepackage{booktabs,arydshln,multirow}
\usepackage{setspace}
\usepackage{xcolor}
\usepackage{stmaryrd}
\usepackage[utf8]{inputenc}
\usepackage{lineno,hyperref}
\usepackage[top=1in, bottom=1in, left=1in, right=1in]{geometry}

\newtheorem{theorem}{Theorem}[section]

\newtheorem{lemma}[theorem]{Lemma}

\newtheorem{definition}[theorem]{Definition}
\newtheorem{remark}[theorem]{Remark}

\modulolinenumbers[10]

\journal{Journal of \LaTeX\ Templates}

\begin{document}

\begin{frontmatter}

\title{Convolution quadrature for Hadamard fractional calculus and correction methods for the subdiffusion with singular source terms}

\author[mymainaddress]{Baoli Yin\corref{mycorrespondingauthor}}
\cortext[mycorrespondingauthor]{Corresponding author}
\ead{baolimath@126.com}

\author[mymainaddress]{Guoyu Zhang}
\ead{guoyu\_zhang@imu.edu.cn}

\author[mymainaddress]{Yang Liu}
\ead{mathliuyang@imu.edu.cn}

\author[mymainaddress]{Hong Li}
\ead{smslh@imu.edu.cn}

\address[mymainaddress]{School of Mathematical Sciences, Inner Mongolia University, Hohhot 010021, China}

\begin{abstract}
The convolution quadrature method originally developed for the Riemann-Liouville fractional calculus is extended in this work to the Hadamard fractional calculus by using the exponential type meshes.
Local truncation error analysis is presented for singular solutions.
By adopting the fractional BDF-$p(1\leq p \leq 6)$ for the Caputo-Hadamard fractional derivative in solving subdiffusion problem with singular source terms, and using the finite element method to discretize the space variable, we carry out the sharp error analysis rigorously and obtain the optimal accuracy by the novel correction technique.
Our correction method is a natural generalization of the one developed for subdiffusion problems with smooth source terms.
Numerical tests confirm the correctness of our theoretical results.
\end{abstract}


\end{frontmatter}


\section{Introduction}
The logarithmic nature in some process is observed in mechanics and engineering.
Unlike the anomalous diffusion which is characterized by a power-law growth of the mean square displacement (MSD) in time \cite{metzler2000random}, the strongly anomalous diffusion has a super-heavy tail in the mean waiting time with the MSD growing logarithmically \cite{yang2022well}. 
This nature can lead to nonlocal operators featured by the logarithmic kernel, e.m., the well-known Hadamard fractional derivative (proposed in 1892 \cite{hadamard1892essai}) or the Caputo-Hadamard (CH) fractional derivative \cite{jarad2012caputo}, which have been successfully applied in  ultra slow kinetics such as in rheology \cite{garra2017generalization}.
More application of Hadamard type calculus in biology, elasticity or turbulent flow can be found in \cite{awadalla2019new,ahmad2017hadamard,wang2019explicit} and references therein.
\par
The mathematical exploration of Hadamard fractional calculus is abundant in literature.
Kilbas and Butzer studied the semigroup properties and the Mellin transform for Hadamard fractional calculus in \cite{kilbas2001hadamard,butzer2002fractional}.
By using the fixed point theory, the authors in \cite{abbas2017survey} carried out the proof of existence of weak solutions of some Hadamard fractional differential systems.
Li and Li in \cite{ma2017hadamard} developed further some reciprocal properties of the Hadamard fractional operator, and proposed the definite conditions of some Hadamard fractional equations.
In \cite{ma2018finite}, Ma and Li revealed the relations between Hadamard fractional derivatives and the finite part integrals in Hadamard sense.
The asymptotic properties of solutions of nonlinear CH
fractional differential equations were explored in \cite{graef2017asymptotic}.
Based on the modified Laplace transform and the finite Fourier sine transform, the authors in \cite{li2020mathematical} obtained the analytic solution to the linear CH fractional models and studied the regularity and logarithmic decay of the solution.
For other aspects about the blow-up phenomenon, the comparison principles and eigenfunctions of Hadamard fractional differential operators, readers can refer \cite{ma2019blow,ma2019comparison,ma2020kinetics}.
\par
In contrast to the classical local operator, the Hadamard fractional calculus operators are featured by the nonlocality and singularity at initial time, both of which render numerically solving problems with such operators extremely challenging.
In history, Gohar et al. \cite{gohar2020caputo,gohar2020finite} developed Euler/$L_{\log,1}$ and predictor-corrector methods for CH fractional differential equations.
Green et al. \cite{green2021numerical} considered predictor-corrector numerical methods on nonuniform meshes under the weak smoothness assumptions of the CH fractional derivative.
In \cite{green2022detailed}, Green and Yan extended the classical Adams methods to models with CH fractional derivatives and established the detailed error analysis.
Li et al. \cite{li2020mathematical} proposed the fully discrete scheme by adopting a difference formula on nonuniform meshes for the time CH fractional derivative in combination with the local discontinuous Galerkin method in spatial direction.
The stability and error estimates were carried out in detail.
For the CH fractional sub-diffusion equations, Wang et al. \cite{wang2022second} built a second-order scheme on nonuniform meshes by resorting to the $L_{\log,2-1_{\sigma}}$ interpolation approximation to the CH fractional derivative and developed sharp error analysis.
Particularly, the authors in \cite{fan2022numerical} proposed three high-order numerical formulas known as L1-2,  L2-1${}_{\sigma}$ and H2N2 formula for the CH fractional derivatives on uniform meshes and uniform partition in the logarithmic sense, respectively.
It is worth mentioning that the mesh which is a uniform partition in the logarithmic sense is also known as exponential type nonuniform mesh \cite{ou2022mathematical} defined by
\begin{equation*}
a=t_0<t_1<t_2<\cdots<t_n<\cdots<t_N=T,\quad
t_n=a\bigg(\frac{T}{a}\bigg)^{\frac{n}{N}}.
\end{equation*}
Such kind of exponential type nonuniform mesh can also be obtained by using a logarithmic transformation \cite{zheng2021logarithmic}.
\par
Undoubtedly, the numerical formulas mentioned above for Hadamard fractional derivatives were inspired to some extent by their counterparts for Caputo/Riemann-Liouville fractional derivatives.
The well-known two groups of difference formulas, i.e. the L-type formulas (L1, L1-2, L2-1${}_{\sigma}$, etc., see \cite{sun2006fully,lin2007finite,alikhanov2015new,liao2020second,feng2019finite} for example) and convolution quadrature (CQ) methods (see \cite{lubich1986discretized,ding2021development,zeng2015numerical,zhang2022convergence,liu2021unified,yin2022efficient} for example) play the essential role in designing numerical schemes.
However, we emphasize that, to the best of our knowledge, none high-order CQ methods have been extended to the Hadamard fractional calculus, which will be our main motivation in this work.
It is known that the solution of CH fractional differential equations is singular at initial time, even if the source term is smooth enough.
We shall in this paper deal with this singularity even if the source term is of weak singular, aiming to yield the desired optimal accuracy by resorting to the newly developed correction technique.
It is notable that such novel correction technique can also be applied to solving Caputo fractional differential equations with singular source terms.
Some other technique to surmount such singularity of source terms can be found in \cite{shi2023high,zhou2023crank}.
\par
The contribution in this paper is from three aspects:
\begin{itemize}
\item[(i)] The high-order (up to 6th) CQ methods is extended to approximate Hadamard fractional calculus operators.
\item[(ii)] Novel correction scheme is developed to maintain the optimal accuracy for the Caputo-Hadamard fractional subdiffusion problem with singular source terms.
\item[(iii)] Sharp error analysis are carried out in detail with numerical confirmation.
\end{itemize}
\par
The outline of the rest of the paper is as follows.
In Sec. \ref{sec.pre}, some preliminaries to the subject and useful tools are introduced for later theoretical analysis.
In Sec. \ref{sec.dis}, local truncation error of the CQ methods is presented based on the exponential type meshes for the Hadamard fractional calculus operators.
By using the fractional BDF-$p(1\leq p \leq 6)$ to approximate the CH fractional derivative, sharp error analysis of the finite element method for the CH fractional subdiffusion problem with singular source terms is presented in Sec. \ref{sec.app}.
Numerical experiments are implemented in Sec. \ref{sec.num} to verify the results of the theoretical analysis and finally, some comments are given in Sec. \ref{sec.con}.
\par
Through out the article, by $A \lesssim B$ we mean that there exists a positive constant $c$ which is independent of the mesh size $h$ or $\tau$, satisfying $A \leq c B$.
\section{Preliminaries}\label{sec.pre}

Denote by $f(t)$ and $g(t)$ functions defined on $[a,+\infty)$.
To be simplicity, define $\delta_t:=t\frac{\mathrm{d}}{\mathrm{d}t}$ and assume $m\in \mathbb{Z}^+$ satisfying $m-1<\alpha<m$.

\begin{definition}
For a given function $f(t)$ with $t>a>0$, define that
\begin{itemize}
\item[{\rm (i)}] the $\alpha$th-order Hadamard fractional integral \cite{hadamard1892essai,li2020mathematical}: 
\begin{equation*}\begin{split}
({}_HI_a^{\alpha}f)(t)=\frac{1}{\Gamma(\alpha)}\int_a^t \bigg(\log\frac{t}{s}\bigg)^{\alpha-1}f(s)\frac{\mathrm{d}s}{s},
\end{split}\end{equation*}
\item[{\rm (ii)}] the $\alpha$th-order Hadamard fractional derivative \cite{hadamard1892essai,li2020mathematical}:
\begin{equation*}\begin{split}
({}_HD_a^\alpha f)(t)=\frac{\delta_t^m}{\Gamma(m-\alpha)}\int_a^t
\bigg(\log\frac{t}{s}\bigg)^{m-\alpha-1}f(s)\frac{\mathrm{d}s}{s},
\end{split}\end{equation*}
\item[{\rm (iii)}] the $\alpha$th-order CH fractional derivative \cite{jarad2012caputo,li2020mathematical}: 
\begin{equation*}\begin{split}
({}_{CH}D_a^\alpha f)(t)=\frac{1}{\Gamma(m-\alpha)}\int_a^t
\bigg(\log\frac{t}{s}\bigg)^{m-\alpha-1}\delta_s^mf(s)\frac{\mathrm{d}s}{s}.
\end{split}\end{equation*}
\end{itemize}
\end{definition}
\begin{definition}
(See \cite{podlubny1999introduction})
For a given function $f(t)$ with $t>0$, define that
\begin{itemize}
\item[{\rm (i)}] the $\alpha$th-order Riemann-Liouville fractional integral: 
\begin{equation*}\begin{split}
({}_RI_0^\alpha f)(t)=\frac{1}{\Gamma(\alpha)}\int_0^t
(t-s)^{\alpha-1}f(s)\mathrm{d}s,
\end{split}\end{equation*}
\item[{\rm (ii)}] the $\alpha$th-order Riemann-Liouville fractional derivative:
\begin{equation*}\begin{split}
({}_RD_0^\alpha f)(t)=\frac{1}{\Gamma(m-\alpha)}\frac{\mathrm{d}^m}{\mathrm{d}t^m}\int_0^t
(t-s)^{m-\alpha-1}f(s)\mathrm{d}s,
\end{split}\end{equation*}
\item[{\rm (iii)}] the $\alpha$th-order Caputo fractional derivative: 
\begin{equation*}\begin{split}
({}_CD_0^\alpha f)(t)=\frac{1}{\Gamma(m-\alpha)}\int_0^t
(t-s)^{m-\alpha-1}f^{(m)}(s)\mathrm{d}s.
\end{split}\end{equation*}
\end{itemize}
\end{definition}
\begin{lemma}(Theorem 3.4 in \cite{kilbas2001hadamard})
The Hadamard fractional derivative and CH fractional derivative have the following relationship
\begin{equation*}
({}_HD_a^\alpha f)(t)=({}_{CH}D_a^\alpha f)(t)
+\sum_{k=0}^{m-1}\frac{(\delta_t^k f)(a)}{\Gamma(k+1-\alpha)}\bigg(\log\frac{t}{a}\bigg)^{k-\alpha}.
\end{equation*}
Particularly, if $\alpha \in (0,1)$, there holds
\begin{equation}\label{Pre.1}
{}_H D_a^\alpha \big(f-f(a)\big)={}_{CH} D_a^\alpha f.
\end{equation}
\end{lemma}
\par
We conclude this section by introducing some useful tools in analyzing the Hadamard fractional calculus.
\begin{definition}(Definitions 2.1 and 2.2 in \cite{li2020mathematical})
The modified Laplace transform $\mathcal{L}\{\cdot\}$ of $f(t)$ is defined by
\begin{equation*}\begin{split}
\widehat{f}(z)=\mathcal{L}\{f(t)\}(z)
=\int_a^\infty e^{-z\log\frac{t}{a}}f(t)\frac{\mathrm{d}t}{t},\quad z\in \mathbb{C}.
\end{split}\end{equation*}
The corresponding inverse transform is defined by
\begin{equation*}\begin{split}
f(t)=\mathcal{L}^{-1}\{\widehat{f}(z)\}(t)=\frac{1}{2\pi\rm{i}}\int_{\mathcal{C}}e^{z\log\frac{t}{a}}\widehat{f}(z)\mathrm{d}z,\quad \rm{i}^2=-1,
\end{split}\end{equation*}
where the contour $\mathcal{C}$ will be clarified later.
\end{definition}
The modified Laplace transform of the CH fractional derivative can be presented in a similar manner as the standard Laplace transform of the Caputo fractional derivative, which takes the form (see Lemma 2.1 in \cite{li2020mathematical})
\begin{equation}\begin{split}\label{Pre.1.1}
\mathcal{L}\{{}_{CH}D_a^\alpha f\}(z)&=z^\alpha\mathcal{L}\{f\}-\sum_{k=0}^{m-1}z^{\alpha-k-1}(\delta_t^kf)(a).
\end{split}\end{equation}
One can also directly verify the following modified Laplace transform of powers of logarithmic functions 
\begin{equation}\begin{split}\label{Pre.1.1.1}
\mathcal{L}\bigg\{\bigg(\log\frac{t}{a}\bigg)^s\bigg\}(z)&=\frac{\Gamma(s+1)}{z^{s+1}}.
\end{split}\end{equation}
\begin{definition}(Definition 2.3 in \cite{li2020mathematical})
The modified convolution of $f(t)$ and $g(t)$ is defined by
\begin{equation*}\begin{split}
(f\ast g)(t)=\int_a^tf\bigg(a\frac{t}{s}\bigg)g(s)\frac{\mathrm{d}s}{s}.
\end{split}\end{equation*}
\end{definition}
The modified convolution operator $\ast$ is closely related to the classical one $\overline{\ast}$ in the sense that for given $t$, if we define $\overline{t}$ such that $t=ae^{\overline{t}}$ and set $F(\overline{t}):=f(t)$, $G(\overline{t}):=g(t)$, there holds
\begin{equation}\begin{split}\label{Pre.1.2}
[f(t)\ast g(t)](t)=[F(\overline{t}) \overline{\ast} G(\overline{t})](\overline{t}).
\end{split}\end{equation}
The modified Laplace transform and modified convolution meet the  property (Proposition 2.2 in \cite{li2020mathematical}) that
\begin{equation*}\begin{split}
\mathcal{L}\{f\ast g\}=\mathcal{L}\{f\}\mathcal{L}\{g\},\quad
\mathcal{L}^{-1}\{\widehat{f}(z)\widehat{g}(z)\}=(f\ast g)(t).
\end{split}\end{equation*}
The following theorem is a modified version of the Taylor theorem with integral remainder which will be used later in error estimates.
\begin{theorem}\label{thm.0}
Let $f^{(k)}$ be absolutely continuous on $[a,t]$, then
\begin{equation}\begin{split}\label{Pre.2}
f(t)=f(a)+\big(\delta_t f\big)\big|_{t=a}\log\frac{t}{a}
+\frac{1}{2!}\big(\delta_t^2 f\big)\big|_{t=a}\bigg(\log\frac{t}{a}\bigg)^2
+\cdots+\frac{1}{k!}\big(\delta_t^k f\big)\big|_{t=a}\bigg(\log\frac{t}{a}\bigg)^k
+R(t),
\end{split}\end{equation}
where the remainder $R(t)$ is given by
\begin{equation*}\begin{split}
R(t)=\frac{1}{k!}\bigg(\log\frac{t}{a}\bigg)^k\ast \delta_t^{k+1}f
=\frac{1}{k!}\int_{a}^t\bigg(\log\frac{t}{s}\bigg)^k\delta_s^{k+1}f(s)\frac{\mathrm{d}s}{s}.
\end{split}\end{equation*}
\end{theorem}

\section{Discretization on exponential type mesh}\label{sec.dis}
For given final time $T>0$ and the time variable $t \in [a,T]$, by using the transform $\overline{t}=\log\frac{t}{a} \in [0,\log\frac{T}{a}]$, we introduce the uniform mesh for the $\overline{t}$-coordinate system \cite{zheng2021logarithmic} :
\begin{equation}\label{Dis.1}
0=\overline{t}_0<\overline{t}_1<\overline{t}_2<\cdots<\overline{t}_n<\cdots<\overline{t}_N=\overline{T}:=\log\frac{T}{a},\quad
\overline{t}_n=n\overline{\tau},\quad
\overline{\tau}=\frac{1}{N}\log\frac{T}{a},\quad
N \in \mathbb{Z}^+.
\end{equation}
The corresponding exponential type mesh for the $t$-coordinate system is then formulated as
\begin{equation}\label{Dis.2}
a=t_0<t_1<t_2<\cdots<t_n<\cdots<t_N=T,\quad
t_n=ae^{\overline{t}_n}=a\bigg(\frac{T}{a}\bigg)^{\frac{n}{N}}.
\end{equation}
To simplify the notation, we extend the operator ${}_RD_0^\alpha$ to the case $\alpha<0$ standing for ${}_R I_0^{-\alpha}$.
On the uniform mesh (\ref{Dis.1}), in accordance to the CQ theory, one can take $\omega_p(\zeta;\alpha)$ as a generating function for the operator ${}_RD_0^\alpha$  which is of $p$th order convergent (Definition 2.3 in \cite{lubich1986discretized}), i.e., satisfying the following properties (see Theorem 2.5 in  \cite{lubich1986discretized}):
\begin{equation}\label{Dis.3}
\text{(i)}\quad \omega_n^{(\alpha)}=O(n^{-\alpha-1}),\quad\quad
\text{(ii)}\quad \overline{\tau}^{-\alpha}\omega_p(e^{-\overline{\tau}};\alpha)-1=O(\overline{\tau}^p),
\end{equation}
where $\omega_n^{(\alpha)}$ is the $n$th coefficient of the series of $\omega_p(\zeta;\alpha)=\sum_{k=0}^\infty\omega_k^{(\alpha)}\zeta^k$.
For example, two groups of the classical CQ methods include the $p$th order fractional backward difference formulas (fractional BDF-$p$)
\begin{equation}\label{Dis.3.1}
\omega_p(\zeta;\alpha)=\psi_p(\zeta)^\alpha,\quad \psi_p(\zeta)=\sum_{j=1}^p\frac{1}{j}(1-\zeta)^j,
\end{equation}
and the generalized Newton-Gregory formula
\begin{equation*}
\omega_p(\zeta;\alpha)=\psi_p(\zeta)^\alpha,\quad 
\psi_p(\zeta)=(1-\zeta)\big[\gamma_0+\gamma_1(1-\zeta)+\gamma_2(1-\zeta)^2+\cdots+\gamma_{p-1}(1-\zeta)^{p-1}\big]^{\frac{1}{\alpha}},
\end{equation*}
where the coefficients $\gamma_i$ are from the series $\sum_{i=0}^\infty\gamma_i(1-\zeta)^i=\big(\frac{\log \zeta}{\zeta-1}\big)^\alpha$.
More CQ methods can be found in \cite{lubich1986discretized}.
To simplify the notation, define
\begin{equation}\begin{split}
\psi_{p,\overline{\tau}}(\zeta)=
\frac{\psi_{p}(\zeta)}{\overline{\tau}}.
\end{split}\end{equation}

\begin{lemma}(Definition 2.3 and Theorem 2.5 in \cite{lubich1986discretized})
Let $U(\overline{t})=\overline{t}^\sigma$ with $\sigma>0$.
Assume the generating function $\omega_p(\zeta,\alpha)$ satisfies (i) and (ii) in (\ref{Dis.3}) with coefficients $\{\omega_n^{(\alpha)}\}$.
Define the discrete convolution $D_{\overline{\tau},n}^\alpha U=\overline{\tau}^{-\alpha}\sum_{k=0}^n\omega_{n-k}^{(\alpha)}U(\overline{t}_k)$ on the uniform mesh (\ref{Dis.1}).
There holds
\begin{equation}\label{Dis.4}
D_{\overline{\tau},n}^\alpha U-({}_RD_0^\alpha U)(\overline{t}_n)=O(\overline{t}_n^{-\alpha-1}\overline{\tau}^{\sigma+1})+O(\overline{t}_n^{-\alpha+\sigma-p}\overline{\tau}^p).
\end{equation}
\end{lemma}
\begin{theorem}\label{thm.1}
Let $\alpha \in \mathbb{R}$ and $\alpha \neq 0,1,2,\cdots$.
Take $u(t)=\big(\log\frac{t}{a}\big)^\sigma$ with $\sigma>0$.
Assume the generating function $\omega_p(\zeta,\alpha)$ satisfies (i) and (ii) in (\ref{Dis.3}) with coefficients $\{\omega_n^{(\alpha)}\}$.
Define the discrete convolution $\mathbb{D}_{\overline{\tau},n}^\alpha u=\overline{\tau}^{-\alpha}\sum_{k=0}^n\omega_{n-k}^{(\alpha)}u({t}_k)$ on the exponential type mesh (\ref{Dis.2}).
There holds
\begin{equation}\label{Dis.5}
\mathbb{D}_{\overline{\tau},n}^\alpha u-({}_HD_a^\alpha u)(t_n)=
O\bigg(\big(\log(t_n/a)\big)^{-\alpha-1}\overline{\tau}^{\sigma+1}\bigg)
+O\bigg(\big(\log(t_n/a)\big)^{-\alpha+\sigma-p}\overline{\tau}^p\bigg).
\end{equation}
\begin{proof}
Let $\overline{t}=\log\frac{t}{a}$ and $u(t)=u(ae^{\overline{t}})=:U(\overline{t})$.
We first show that $({}_HD_a^\alpha u)(t)=({}_RD_0^\alpha U)(\overline{t})$.
\par
Case I: $\alpha<0$.
\par
In this case, ${}_HD_a^\alpha$ represents the Hadamard fractional integral operator ${}_H I_a^{-\alpha}$.
By definition, one gets
\begin{equation*}\begin{split}
({}_HD_a^\alpha u)(t)&=\frac{1}{\Gamma(-\alpha)}\int_a^t\bigg(\log\frac{t}{s}\bigg)^{-\alpha-1}u(s)\frac{\mathrm{d}s}{s}
\xlongequal{s=ae^w}
\frac{1}{\Gamma(-\alpha)}\int_0^{\log\frac{t}{a}}\bigg(\log\frac{t}{a}-w\bigg)^{-\alpha-1}u(ae^w)\mathrm{d}w
\\
&=\frac{1}{\Gamma(-\alpha)}\int_0^{\overline{t}}\big(\overline{t}-w\big)^{-\alpha-1}U(w)\mathrm{d}w
\\
&=({}_RD_0^\alpha U)(\overline{t}).
\end{split}\end{equation*}
\par
Case II: $m-1<\alpha<m$ for some $m\in \mathbb{Z}^+$.
\par
In this case, by definition we have
\begin{equation}\begin{split}
({}_HD_a^\alpha u)(t)&=\frac{1}{\Gamma(m-\alpha)}\bigg(t\frac{\mathrm{d}}{\mathrm{d}t}\bigg)^m\int_a^t\bigg(\log\frac{t}{s}\bigg)^{m-\alpha-1}u(s)\frac{\mathrm{d}s}{s}
\\
&\xlongequal{s=ae^w}
\frac{1}{\Gamma(m-\alpha)}\bigg(t\frac{\mathrm{d}}{\mathrm{d}t}\bigg)^m\int_0^{\log\frac{t}{a}}\bigg(\log\frac{t}{a}-w\bigg)^{m-\alpha-1}u(ae^w)\mathrm{d}w
\\
&=\frac{1}{\Gamma(m-\alpha)}\bigg(ae^{\overline{t}}\frac{\mathrm{d}\overline{t}}{\mathrm{d}t}\frac{\mathrm{d}}{\mathrm{d}\overline{t}}\bigg)^m
\int_0^{\overline{t}}(\overline{t}-w)^{m-\alpha-1}U(w)\mathrm{d}w
\\
&=\frac{1}{\Gamma(m-\alpha)}\bigg(\frac{\mathrm{d}}{\mathrm{d}\overline{t}}\bigg)^m
\int_0^{\overline{t}}(\overline{t}-w)^{m-\alpha-1}U(w)\mathrm{d}w
\\
&=({}_RD_0^\alpha U)(\overline{t}).
\end{split}\end{equation}
The estimate (\ref{Dis.5}) follows readily from (\ref{Dis.4}) since 
$$\mathbb{D}_{\overline{\tau},n}^\alpha u=\overline{\tau}^{-\alpha}\sum_{k=0}^n\omega_{n-k}^{(\alpha)}u({t}_k)=\overline{\tau}^{-\alpha}\sum_{k=0}^n\omega_{n-k}^{(\alpha)}U(\overline{t}_k)=D_{\overline{\tau},n}^\alpha U,$$
combined with the fact $\overline{t}_n=\log\frac{t_n}{a}$.\end{proof}
\end{theorem}

\section{Application in CH fractional subdiffusion problem with sharp estimates}\label{sec.app}
In this section, we adopt the high-order CQ method (fractional BDF-$p$) to the CH fractional subdiffusion problem with a singular source term, which can be formulated as
\begin{equation}\begin{split}\label{App.1}
\begin{cases}
{}_{CH}D_a^\alpha u(\boldsymbol x,t)-\Delta u(\boldsymbol x,t)=f(\boldsymbol x,t),\quad t \in (a,T],\quad \boldsymbol x\in \Omega, \quad \alpha\in(0,1),
\\
u(\boldsymbol x,a)=v(\boldsymbol x),\quad \boldsymbol x\in\Omega,\\
u(\boldsymbol x,a)=0,\quad \boldsymbol x\in\partial\Omega,\quad t\in(a,T],
\end{cases}
\end{split}\end{equation}
where $\Omega \subset \mathbb{R}^d$ is a bounded convex polygonal domain and $\boldsymbol x=(x_1,\cdots,x_d)^T$ denotes the spatical variable.
The source term $f(\boldsymbol x,t)$ is singular such that
\begin{equation}\begin{split}\label{App.1.1}
f(\boldsymbol x,t)=f(\boldsymbol x,a)+\bigg(\log\frac{t}{a}\bigg)^\beta g(\boldsymbol x,t),\quad \beta \in [0,1),\quad \text{$g$ is smooth enough with respect to $t$}.
\end{split}\end{equation}
$\Delta$ stands for the Laplace operator defined on $D(\Delta):=H_0^1(\Omega) \cap H^2(\Omega)$.
The initial conditions $v$ is smooth enough in the sense $v \in D(\Delta)$.
See Remark \ref{rem.3} for the discussion of nonsmooth initial conditions.
\par
Our main interest lies in analyzing time discretization using time-stepping methods. Therefore, we begin by formulating the space semidiscrete scheme using the finite element method. Let $\mathcal{T}_h$ be a shape regular, quasi-uniform triangulation of the domain $\Omega$ into  $d$-simplexes where $h$ denotes the mesh size. We approximate the solution in each simplex $e$ by a $k$th order polynomial function.
The finite element space is defined by
\begin{equation*}\label{App.2}\begin{split}
S_h=\{\chi_h \in H_0^1(\Omega): \chi_h |_e \in \mathcal{P}_k(e),~ e \in \mathcal{T}_h\},
\end{split}\end{equation*}
where $\mathcal{P}_k(e)$ denotes the space of all polynomials on $e$ with order at most $k$.

Introduce the $L^2(\Omega)$ projection $P_h:L^2(\Omega)\to S_h$, the Ritz projection $R_h:H_0^1(\Omega)\to S_h$ and the discrete Laplacian $\Delta_h: S_h\to S_h$, respectively, as follows,
\begin{equation}\label{App.3}\begin{split}
(P_h \phi,\chi_h)&=(\phi,\chi_h), \quad\forall \phi \in L^2(\Omega),~ \forall \chi_h \in S_h,
\\
(\nabla R_h \phi,\nabla \chi_h)&=(\nabla \phi,\nabla \chi_h),\quad \forall \phi \in H_0^1(\Omega),~ \forall \chi_h \in S_h,
\\
(\Delta_h \phi_h,\chi_h)&=-(\nabla \phi_h, \nabla \chi_h), \quad \forall \phi_h, \chi_h \in S_h.
\end{split}\end{equation}
We mention that as a sectorial operator, the discrete Laplacian $\Delta_h$ meets the resolvent estimate
\begin{equation}\begin{split}\label{App.3.1}
\|(z-\Delta_h)^{-1}\|\lesssim |z|^{-1},\quad \forall z \in \Sigma_\theta,\quad \text{for some } \theta >\frac{\pi}{2},
\end{split}\end{equation}
where $\Sigma_\theta$ denotes the open sector $\{z\in \mathbb{C}:|\arg z|<\theta,z\neq 0\}$.
Therefore, one readily gets $\|(z^\alpha-\Delta_h)^{-1}\|\lesssim |z|^{-\alpha}$ provided $\alpha \in (0,1)$ and $z \in \Sigma_\theta$.
\subsection{Fully discrete scheme}
The space semi-discrete scheme of problem (\ref{App.1}) is to find $u_h \in S_h$ such that
\begin{equation}\begin{split}\label{App.4}
{}_{CH}D_a^\alpha u_h(\boldsymbol x,t)-\Delta_h u_h(\boldsymbol x,t)
=P_hf(\boldsymbol x,t)
=P_h f(\boldsymbol x,a)+\bigg(\log\frac{t}{a}\bigg)^\beta P_h g(\boldsymbol x,t),
\end{split}\end{equation}
with the initial condition $u_h(\boldsymbol x,a)=v_h(\boldsymbol x):=R_h v(\boldsymbol x)$.
Let $f_h(\boldsymbol x,t):=P_h f(\boldsymbol x,t)$, $g_h(\boldsymbol x,t):=P_h g(\boldsymbol x,t)$ and $w_h(\boldsymbol x,t):=u_h(\boldsymbol x,t)-v_h(\boldsymbol x)$.
By the relation (\ref{Pre.1}), we can rewrite the problem (\ref{App.4}) into the form
\begin{equation}\begin{split}\label{App.5}
{}_{H}D_a^\alpha w_h(\boldsymbol x,t)-\Delta_h w_h(\boldsymbol x,t)
=\Delta_h v_h(\boldsymbol x)
+f_h(\boldsymbol x,t).
\end{split}\end{equation}
To simplify the notation, for a general function $\psi(\boldsymbol x,t)$, let $\psi^n(\boldsymbol x):=\psi(\boldsymbol x,t_n)$ and the space variable will be omitted.
\par
By adopting the fractional BDF-$p$ (\ref{Dis.3.1}) to discretize the time variable on the exponential type mesh (\ref{Dis.2}), we propose the fully discrete scheme with correction terms: 
\begin{equation}\begin{split}\label{App.5.1}
\mathbb{D}_{\overline{\tau},n}^\alpha W_h-\Delta_h W_h^n
&=
\Delta_hv_h+f_h^n
+\underbrace{
b_n^{(p)} (\Delta_h v_h+f_h^0)
+\sum_{j=0}^{p-2} d_{j,n}^{(p,\beta)}\overline{\tau}^{j+\beta} \delta_s^j g_h(a)}_{\text{Correction terms}}
,\quad 1\leq n \leq p-1,
\\
\mathbb{D}_{\overline{\tau},n}^\alpha W_h-\Delta_h W_h^n
&=
\Delta_h v_h+f_h^n,\quad p\leq n \leq N,
\end{split}\end{equation}
where $\delta_s^j g_h(a)=\big(\delta_s^j g_h(s)\big)\big|_{s=a}$ and the coefficients $b_n^{(p)}$  and $d_{j,n}^{(p,\beta)}$ will be determined later on.
We adopt the convention that if the upper index of the summation is less than the lower one, the summation is nothing but zero.
\par
It is notable that BDF-$p$ $(1\leq p \leq 6)$ are A-$(\vartheta_p)$ stable \cite{wanner1996solving} with $\vartheta_p=90^{\circ}, 90^{\circ}, 86.03^{\circ}, 73.35^{\circ}, 51.84^{\circ}, 17.84^{\circ}$ for $p=1,2,3,4,5,6$, respectively.
Denote by $\mathbb{U}_r=\{z\in \mathbb{C}:|z|<r\}$.
The function $\psi_p(\zeta)$ has no zeros except $\zeta=1$ in $\mathbb{U}_r$ for some $r>1$ \cite{creedon1975stability}, which leads to the invertibility of the operator 
$\overline{\tau}^{-\alpha}\psi_p(\zeta)^\alpha-\Delta_h$ for $\zeta \in \mathbb{U}_r \setminus \{1\}$ with some $r>1$, satisfying
\begin{equation}\begin{split}\label{App.5.1.1}
\|(\overline{\tau}^{-\alpha}\psi_p(\zeta)^\alpha-\Delta_h)^{-1}\|\lesssim \overline{\tau}^\alpha|\psi_p(\zeta)|^{-\alpha}.
\end{split}\end{equation}

\subsection{Calculation for the coefficients $b_n^{(p)}$ and $d_{j,n}^{(p,\beta)}$}
In accordance with the error estimates in the next subsection, the coefficients $b_n^{(p)}$ and $d_{j,n}^{(p,\beta)}$ are determined by justifying the following inequalities:
\begin{equation}\begin{split}\label{App.5.2}
|\mu(\zeta)-1|
 &\lesssim |1-\zeta|^p,
\\
\bigg|
\beta_j(\zeta)
-\frac{\Gamma(j+1+\beta)}{j! \psi_p(\zeta)^{j+1+\beta}}\bigg|&\lesssim 
|1-\zeta|^{p-j-1-\beta},\quad j=0,1,2,\cdots, p-2,
\end{split}\end{equation}
with $\mu(\zeta)$ and $\beta_j(\zeta)$ defined by
\begin{equation*}\begin{split}
\mu(\zeta)&=\psi_p(\zeta)\bigg(\sum_{n=1}^{p-1}b_n^{(p)}\zeta^n
+\frac{\zeta}{1-\zeta}\bigg),
\\
\beta_j(\zeta)&=\sum_{n=1}^{p-1}d_{j,n}^{(p,\beta)}\zeta^n+\frac{1}{j!}{\rm Li}_{-(j+\beta)}(\zeta),
\end{split}\end{equation*}
where  ${\rm Li}_s(\zeta)$ stands for the polylogarithm function
$
{\rm Li}_s(\zeta)=\displaystyle\sum_{j=1}^\infty \frac{\zeta^j}{j^s}
$.
It is notable that ${\rm Li}_s(\zeta)$ can be expanded at $\zeta=1$ into the series (see Theorem 1 in \cite{flajolet1999singularity})
\begin{equation}\begin{split}\label{App.5.3}
{\rm Li}_s(\zeta)=
\Gamma(1-s)\big(-\log \zeta\big)^{s-1}
+\sum_{k=0}^\infty \frac{\mathcal{Z}(s-k)}{k!}(\log \zeta)^k,
\end{split}\end{equation}
where $\mathcal{Z}(s):=\displaystyle\sum_{n=1}^\infty \frac{1}{n^s}$ is the Riemann zeta function which is convergent for $\Re(s)>1$ and can be analytically continued to all $s \in \mathbb{C}\setminus \{1\}$.

The coefficients $b_n^{(p)}$ can be found in Table 1 in \cite{jin2017correction}.
We mention that $b_n^{(p)}$ is exactly $d_{j,n}^{(p,\beta)}$ with $j=\beta=0$. See Table \ref{tab1}.
We next present the derivation of the coefficients $d_{j,n}^{(p,\beta)}$.
Note that the generating function $\psi_p(\zeta)$ for the BDF-$p$ satisfies the following conditions
\begin{equation*}\begin{split}
\psi_p(\zeta)=(1-\zeta)\Theta_p(\zeta),\quad
\Theta_p(\zeta)=\frac{-\log\zeta}{1-\zeta}+O(|1-\zeta|^p),
\end{split}\end{equation*}
which combined with the expension (\ref{App.5.3}),  leads to
\begin{equation}\begin{split}\label{App.5.4}
&\quad\frac{1}{j!}{\rm Li}_{-(j+\beta)}(\zeta)
-\frac{\Gamma(j+1+\beta)}{j! \psi_p(\zeta)^{j+1+\beta}}
\\
&=\frac{1}{j!}\Gamma(1+j+\beta)(1-\zeta)^{-j-1-\beta}
\bigg[\bigg(\frac{-\log\zeta}{1-\zeta}\bigg)^{-j-1-\beta}-\Theta_p(\zeta)^{-j-1-\beta}\bigg]
\\&\quad+
\frac{1}{j!}\sum_{k=0}^\infty\frac{1}{k!}\mathcal{Z}(-j-k-\beta)(\log \zeta)^k
\\
&=
\frac{1}{j!}\sum_{k=0}^\infty\frac{1}{k!}\mathcal{Z}(-j-k-\beta)(\log \zeta)^k
+O(|1-\zeta|^{p-j-1-\beta}).
\end{split}\end{equation}
Assume the first part of the above last step fulfills the expansion at $\zeta=1$ as follows
\begin{equation}\begin{split}
\frac{1}{j!}\sum_{k=0}^\infty\frac{1}{k!}\mathcal{Z}(-j-k-\beta)(\log \zeta)^k
&=\sum_{k=0}^{\infty}c_{j,k}^{(\beta)}(1-\zeta)^k,
\end{split}\end{equation}
where, by using the expansion of $\log \zeta$ at $\zeta=1$, the first five $c_{j,k}^{(\beta)}$ can be calculated by the formulas
\begin{equation*}\begin{split}
c_{j,0}^{(\beta)} &=\frac{1}{j!}\mathcal{Z}(-j-\beta),
\\
c_{j,1}^{(\beta)}&=-\frac{1}{j!}\mathcal{Z}(-j-\beta-1),
\\
c_{j,2}^{(\beta)}&=-\frac{1}{2 j!}\big[\mathcal{Z}(-j-\beta-1)-\mathcal{Z}(-j-\beta-2)\big],
\\
c_{j,3}^{(\beta)}&=-\frac{1}{6 j!}\big[2 \mathcal{Z}(-j-\beta-1)-3 \mathcal{Z}(-j-\beta-2)+\mathcal{Z}(-j-\beta-3)\big],
\\
c_{j,4}^{(\beta)}&=-\frac{1}{24 j!}\big[6 \mathcal{Z}(-j-\beta-1)-11 \mathcal{Z}(-j-\beta-2)+6\mathcal{Z}(-j-\beta-3)-\mathcal{Z}(-j-\beta-4)\big].
\end{split}\end{equation*}
Introduce the coefficients $\eta_{j,n}^{(p)}$ $(0 \leq n \leq p-2)$ such that
\begin{equation*}\begin{split}
\sum_{n=1}^{p-1}d_{j,n}^{(p,\beta)}\zeta^n
&=
\zeta\sum_{n=0}^{p-2} \eta_{j,n}^{(p)}(1-\zeta)^n
=\sum_{n=0}^{p-2} \eta_{j,n}^{(p)}(1-\zeta)^n
-\sum_{n=0}^{p-2} \eta_{j,n}^{(p)}(1-\zeta)^{n+1}
\\&=
\eta_{j,0}^{(p)}
+\sum_{n=1}^{p-2}\big( \eta_{j,n}^{(p)}- \eta_{j,n-1}^{(p)}\big)(1-\zeta)^n
- \eta_{j,p-2}^{(p)}(1-\zeta)^{p-1},
\end{split}\end{equation*}
which, therefore, leads to 
\begin{equation*}\begin{split}
\eta_{j,0}^{(p)}&=-c_{j,0}^{(\beta)},
\\
\eta_{j,n}^{(p)}&=\eta_{j,n-1}^{(p)}-c_{j,n}^{(\beta)},\quad n=1,2,\cdots, p-j-2,
\\
\eta_{j,n}^{(p)}&=0,\quad n=p-j-1,\cdots, p-2.
\end{split}\end{equation*}
and finally results in
\begin{equation*}\begin{split}
d_{j,k+1}^{(p,\beta)}=
(-1)^{k+1}\sum_{n=k}^{p-j-2}
\binom{n}{k}
\sum_{\ell=0}^n c_{j,\ell}^{(\beta)},\quad
j,k=0,1,\cdots,p-2.
\end{split}\end{equation*}

\begin{table}[]
\centering
\caption{Coefficients $d_{j,n}^{(p,\beta)}$ for $\beta=0$ (smooth source term)}\label{tab1}
{\renewcommand{\arraystretch}{1.5}
\begin{tabular}{|cc|ccccc|}
\hline
BDF-$p$                &       & $d_{j,1}^{(p,0)}$     & $d_{j,2}^{(p,0)}$     & $d_{j,3}^{(p,0)}$   & $d_{j,4}^{(p,0)}$     & $d_{j,5}^{(p,0)}$  \\ \hline
$p=2$                  & $j=0$ & $\frac{1}{2}$       &                     &                   &                     &                  \\ \hline
\multirow{2}{*}{$p=3$} & $j=0$ & $\frac{11}{12}$     & $-\frac{5}{12}$     &                   &                     &                  \\
                       & $j=1$ & $\frac{1}{12}$      & 0                   &                   &                     &                  \\ \hline
\multirow{3}{*}{$p=4$} & $j=0$ & $\frac{31}{24}$     & $-\frac{7}{6}$      & $\frac{3}{8}$     &                     &                  \\
                       & $j=1$ & $\frac{1}{6}$       & $-\frac{1}{12}$     & 0                 &                     &                  \\
                       & $j=2$ & 0                   & 0                   & 0                 &                     &                  \\ \hline
\multirow{4}{*}{$p=5$} & $j=0$ & $\frac{1181}{720}$  & $-\frac{177}{80}$   & $\frac{341}{240}$ & $-\frac{251}{720}$  &                  \\
                       & $j=1$ & $\frac{59}{240}$    & $-\frac{29}{120}$   & $\frac{19}{240}$  & 0                   &                  \\
                       & $j=2$ & $\frac{1}{240}$     & $-\frac{1}{240}$    & 0                 & 0                   &                  \\
                       & $j=3$ & $-\frac{1}{720}$     & 0                   & 0                 & 0                   &                  \\ \hline
\multirow{5}{*}{$p=6$} & $j=0$ & $\frac{2837}{1440}$ & $-\frac{2543}{720}$ & $\frac{17}{5}$    & $-\frac{1201}{720}$ & $\frac{95}{288}$ \\
                       & $j=1$ & $\frac{77}{240}$    & $-\frac{7}{15}$     & $\frac{73}{240}$  & $-\frac{3}{40}$     & 0                \\
                       & $j=2$ & $\frac{1}{96}$      & $-\frac{1}{60}$     & $\frac{1}{160}$   & 0                   & 0                \\
                       & $j=3$ & $-\frac{1}{360}$    & $\frac{1}{720}$     & 0                 & 0                   & 0                \\
                       & $j=4$ & 0                   & 0                   & 0                 & 0                   & 0                \\ \hline
\end{tabular}}
\end{table}

\subsection{Sharp error estimates}

\noindent\textbf{Step 1: Representation of the solution of equation (\ref{App.5})}
\par
Using the modified Taylor formula (\ref{Pre.2}) for $g_h(t)$, the source term can be formulated as
\begin{equation}\begin{split}\label{App.6}
f_h(t)=
f_h(a)+
\sum_{j=0}^{p-2}\frac{1}{j!}\delta_t^j g_h(a)\bigg(\log\frac{t}{a}\bigg)^{j+\beta}
+  R_p(t),
\end{split}\end{equation}
where the local truncation error $R_p(t)$ takes the form 
\begin{equation}\begin{split}\label{App.6.1}
R_p(t)&=R_p^{(1)}(t)+R_p^{(2)}(t)
\\&=\frac{1}{(p-1)!}\delta_t^{p-1} g_h(a)\bigg(\log\frac{t}{a}\bigg)^{p-1+\beta}
+\frac{1}{(p-1)!}\bigg(\log\frac{t}{a}\bigg)^\beta \bigg[\bigg(\log\frac{t}{a}\bigg)^{p-1}\ast \delta_t^{p}g_h\bigg].
\end{split}\end{equation}

\par
To analytically represent the solution of (\ref{App.5}),  we combine (\ref{App.5}) with (\ref{App.6}) before taking the modified Laplace transform for both hand-sides of (\ref{App.5}) and using the facts in (\ref{Pre.1.1}) and (\ref{Pre.1.1.1}) to obtain
\begin{equation*}\begin{split}\label{App.7}
z^\alpha \widehat{w}_h-\Delta_h \widehat{w}_h
=z^{-1}(\Delta_h v_h+f_h(a))
+\sum_{j=0}^{p-2}\frac{\Gamma(j+1+\beta)}{j! z^{j+1+\beta}}\delta_t^j g_h(a)+\widehat{R}_p(z).
\end{split}\end{equation*}
Using the modified inverse Laplace transform, the solution $w_h(t)$ can be represented by
\begin{equation}\begin{split}\label{App.7.1}
w_h(t)&=\frac{1}{2\pi\rm{i}}\int_{\Gamma_{\theta,\rho}}e^{z\log\frac{t}{a}}K(z)(\Delta_h v_h+f_h(a))\mathrm{d}z
\\
&\quad+\frac{1}{2\pi\rm{i}}\int_{\Gamma_{\theta,\rho}}e^{z\log\frac{t}{a}}zK(z)\bigg[\sum_{j=0}^{p-2}\frac{\Gamma(j+1+\beta)}{j! z^{j+1+\beta}}\delta_t^j g_h(a)+\widehat{R}_p(z)\bigg]\mathrm{d}z,
\end{split}\end{equation}
where $K(z)=z^{-1}(z^\alpha-\Delta_h)^{-1}$ denotes the kernel function and the contour $\Gamma_{\theta,\rho}$ is defined by
\begin{equation*}\begin{split}
\Gamma_{\theta,\rho}:=\{z\in \mathbb{C}:|z|=\rho, |\arg z|\leq \theta\}\cup \{z\in\mathbb{C}: z=re^{\rm{i}\theta}, r\geq \rho\},
\end{split}\end{equation*}
with the orientation that it is traced out with an increasing imaginary part as one moves along the contour.
In the discrete context, $\Gamma_{\theta,\rho}$ will be truncated with the same orientation such that
\begin{equation}\begin{split}
\Gamma_{\theta,\rho}^{\overline{\tau}}:=\{z=x+{\rm i}y\in \Gamma_{\theta,\rho}:|y|\leq \pi/\overline{\tau}\}.
\end{split}\end{equation}
\\\\
\textbf{Step 2: Representation of the solution of scheme (\ref{App.5.1})}
\par
Multiplying both hand-side of (\ref{App.5.1}) by $\zeta^n$ and summing the index $n$ from $1$ to $\infty$, we obtain
\begin{equation}\begin{split}\label{App.8}
\sum_{n=1}^\infty\zeta^n\mathbb{D}_{\overline{\tau},n}^\alpha W_h
-\sum_{n=1}^\infty\zeta^n\Delta_h W_h^n
&=(\Delta_h v_h+f_h^0)\sum_{n=1}^{p-1}b_n^{(p)}\zeta^n
+\sum_{n=1}^\infty\zeta^n\Delta_h v_h
\\
&\quad+\sum_{n=1}^\infty\zeta^nf_h^n
+\sum_{j=0}^{p-2}\overline{\tau}^{j+\beta} \delta_t^j g_h(a)\sum_{n=1}^{p-1}d_{j,n}^{(p,\beta)}\zeta^n.
\end{split}\end{equation}
Replace $f_h(t_n)$ in (\ref{App.8}) by (\ref{App.6}) and collect terms to yield
\begin{equation}\begin{split}\label{App.9}
\underbrace{
\sum_{n=1}^\infty\zeta^n\mathbb{D}_{\overline{\tau},n}^\alpha W_h
-\sum_{n=1}^\infty\zeta^n\Delta_h W_h^n
}_{I_1}
&=(\Delta_h v_h+f_h^0)\bigg(\sum_{n=1}^{p-1}b_n^{(p)}\zeta^n
+\underbrace{\sum_{n=1}^\infty\zeta^n}_{I_2}\bigg)
+\underbrace{\sum_{n=1}^\infty \zeta^n R_p^n}_{I_3}
\\
&
\quad+\sum_{j=0}^{p-2} \delta_t^j g_h(a)\bigg(\overline{\tau}^{j+\beta}
\sum_{n=1}^{p-1}d_{j,n}^{(p,\beta)}\zeta^n
+\underbrace{\frac{1}{j!} \sum_{n=1}^\infty \bigg(\log\frac{t_n}{a}\bigg)^{j+\beta}\zeta^n}_{I_4}\bigg).
\end{split}\end{equation}
Since $W_h^0=0$ and $R_p^0=0$, there hold
\begin{equation*}\begin{split}\label{App.10}
I_1=\sum_{n=0}^\infty\zeta^n\mathbb{D}_{\overline{\tau},n}^\alpha W_h
-\sum_{n=0}^\infty\zeta^n\Delta_h W_h^n
=\overline{\tau}^{-\alpha}\omega(\zeta;\alpha)W_h(\zeta)-\Delta_h W_h(\zeta)
=\big(\overline{\tau}^{-\alpha}\psi_p(\zeta)^\alpha-\Delta_h\big)W_h(\zeta)
\end{split}\end{equation*}
and
\begin{equation}\begin{split}\label{App.11}
I_3=\sum_{n=0}^\infty \zeta^n R_p^n=:\widetilde{R}_p(\zeta).
\end{split}\end{equation}
For the term $I_2$, one readily gets
\begin{equation*}\begin{split}\label{App.12}
I_2=-1+\sum_{n=0}^\infty\zeta^n=\frac{\zeta}{1-\zeta}.
\end{split}\end{equation*}
Thanks to $t_n=ae^{\overline{t}_n}=ae^{n\overline{\tau}}$, the term $I_4$ can be formulated as
\begin{equation}\begin{split}\label{App.13}
I_4=\frac{\overline{\tau}^{j+\beta}}{j!} \sum_{n=1}^\infty n^{j+\beta} \zeta^n
=\frac{\overline{\tau}^{j+\beta}}{j!} {\rm Li}_{-(j+\beta)}(\zeta).
\end{split}\end{equation}
Combining (\ref{App.9})-(\ref{App.13}) and using the kernel function $K(z)=z^{-1}(z^\alpha-\Delta_h)^{-1}$, we have
\begin{equation}\begin{split}\label{App.13.1}
W_h(\zeta)&=
\overline{\tau}^{-1}K(\psi_{p,\overline{\tau}}(\zeta)) \mu(\zeta)(\Delta_h v_h+f_h^0)
\\
&\quad+\psi_{p,\overline{\tau}}(\zeta) K(\psi_{p,\overline{\tau}}(\zeta))
\bigg(\widetilde{R}_p(\zeta)+\sum_{j=0}^{p-2}\overline{\tau}^{j+\beta} \delta_t^j g_h(a)
\beta_j(\zeta)\bigg).
\end{split}\end{equation}
Then, by the Cauchy's integral formula, we get
\begin{equation}\begin{split}\label{App.13.2}
W_h(\zeta)=\frac{1}{2\pi\rm{i}}\int_{|\zeta|=\varepsilon}
\frac{W_h(\zeta)}{\zeta^{n+1}}\mathrm{d}\zeta,
\end{split}\end{equation}
where $\varepsilon$ is a small enough positive number.
Define a contour $\Gamma_\varepsilon^{\overline{\tau}}:=\{z=-\overline{\tau}^{-1}\log\varepsilon+{\rm i}y: y\in\mathbb{R}, |y|\leq \pi/\overline{\tau}\}$ with an increasing imaginary part.
Setting $\zeta=e^{-z\overline{\tau}}$, $z\in \Gamma_\varepsilon^{\overline{\tau}}$, in combination (\ref{App.13.1}) with (\ref{App.13.2}), one immediately gets
\begin{equation}\begin{split}\label{App.14}
W_h^n&=\frac{\overline{\tau}}{2\pi\rm{i}}\int_{\Gamma_\varepsilon^{\overline{\tau}}}
e^{z\overline{t}_n}W_h(e^{-z\overline{\tau}})\mathrm{d}z
\\
&=\frac{\overline{\tau}}{2\pi\rm{i}}\int_{\Gamma_{\theta,\rho}^{\overline{\tau}}}
e^{z\overline{t}_n}W_h(e^{-z\overline{\tau}})\mathrm{d}z
\\
&=\frac{1}{2\pi\rm{i}}\int_{\Gamma_{\theta,\rho}^{\overline{\tau}}}
e^{z\overline{t}_n}K(\psi_{p,\overline{\tau}}(e^{-z\overline{\tau}})) \mu(e^{-z\overline{\tau}})(\Delta_h v_h+f_h^0)\mathrm{d}z
\\&\quad
+\frac{1}{2\pi\rm{i}}\int_{\Gamma_{\theta,\rho}^{\overline{\tau}}}
e^{z\overline{t}_n}\psi_{p,\overline{\tau}}(e^{-z\overline{\tau}}) K(\psi_{p,\overline{\tau}}(e^{-z\overline{\tau}}))
\bigg[\overline{\tau}\widetilde{R}_p(e^{-z\overline{\tau}})
+\sum_{j=0}^{p-2}\overline{\tau}^{j+1+\beta} \delta_t^j g_h(a)
\beta_j(e^{-z\overline{\tau}})\bigg]\mathrm{d}z.
\end{split}\end{equation}
\\\\
\textbf{Step 3: Some technical lemmas}
\begin{lemma}\label{lem.0}(Lemma B.1 in \cite{jin2017correction})
Let $\alpha \in (0,1)$.
There hold the following estimates for $z \in \Gamma_{\theta,\rho}^{\overline{\tau}}$ that
\begin{equation*}\begin{split}
|z|\lesssim |\psi_{p,\overline{\tau}}(e^{-z\overline{\tau}})|\lesssim |z|,\quad
|\psi_{p,\overline{\tau}}(e^{-z\overline{\tau}})-z|\lesssim \overline{\tau}^p|z|^{p+1},\quad
|\psi_{p,\overline{\tau}}(e^{-z\overline{\tau}})^\alpha-z^\alpha|\lesssim \overline{\tau}^p |z|^{p+\alpha}.
\end{split}\end{equation*}
\end{lemma}
\begin{remark}\label{rem.1}
It is notable that the second and third estimates in Lemma \ref{lem.0} can be generalized to the case that for any $s \in \mathbb{R}$,
$|\psi_{p,\overline{\tau}}(e^{-z\overline{\tau}})^s-z^s|\lesssim \overline{\tau}^p |z|^{p+s}$.
Indeed, assume $s \in (1,2]$, then
\begin{equation}\begin{split}\label{App.14.0.0}
|\psi_{p,\overline{\tau}}(e^{-z\overline{\tau}})^s-z^s|
\lesssim
|\psi_{p,\overline{\tau}}(e^{-z\overline{\tau}})|
|\psi_{p,\overline{\tau}}(e^{-z\overline{\tau}})^{s-1}-z^{s-1}|
+|z|^{s-1}|\psi_{p,\overline{\tau}}(e^{-z\overline{\tau}})-z|
\lesssim
\overline{\tau}^p |z|^{p+s}.
\end{split}\end{equation}
Therefore, the mathematical induction method guarantees (\ref{App.14.0.0}) holds also for any $s>0$.
The case $s=0$ is trivial.
For the case $s<0$, one readily gets
\begin{equation*}\begin{split}
|\psi_{p,\overline{\tau}}(e^{-z\overline{\tau}})^s-z^s|
=|\psi_{p,\overline{\tau}}(e^{-z\overline{\tau}})|^{s}|z|^{s}
|\psi_{p,\overline{\tau}}(e^{-z\overline{\tau}})^{-s}-z^{-s}|
\lesssim
\overline{\tau}^p |z|^{p+s}.
\end{split}\end{equation*}
\end{remark}
\begin{lemma}\label{lem.0.1}
Assume $H_{\overline{\tau}}(z)$ is an approximation to a given function $H(z)$ with the estimate
\begin{equation*}\begin{split}
|H(z)-H_{\overline{\tau}}(z)| \lesssim \overline{\tau}^{s_1}|z|^{s_2},\quad \text{for } z \in \Gamma_{\theta,\rho}^{\overline{\tau}}.
\end{split}\end{equation*}
Moreover, assume $|H(z)| \lesssim |z|^{s_3}$.
$s_i (i=1,2,3)$ denotes some constant.
It holds, for any $z \in \Gamma_{\theta,\rho}^{\overline{\tau}}$ that
\begin{equation}\begin{split}\label{App.14.0.1}
\big\|
zK(z)H(z)-
\psi_{p,\overline{\tau}}(e^{-z\overline{\tau}}) K(\psi_{p,\overline{\tau}}(e^{-z\overline{\tau}}))H_{\overline{\tau}}(z)
\big\|
\lesssim
\overline{\tau}^p |z|^{p-\alpha+s_3}
+\overline{\tau}^{s_1}|z|^{s_2-\alpha}.
\end{split}\end{equation}
\begin{proof}
By the fact $zK(z)=(z^\alpha-\Delta_h)^{-1}$ and that
\begin{equation*}\begin{split}
(z^\alpha-\Delta_h)^{-1}
-\big(\psi_{p,\overline{\tau}}(e^{-z\overline{\tau}})-\Delta_h\big)^{-1}
=
(z^\alpha-\Delta_h)^{-1}
\big(\psi_{p,\overline{\tau}}(e^{-z\overline{\tau}})-z^\alpha\big)
\big(\psi_{p,\overline{\tau}}(e^{-z\overline{\tau}})-\Delta_h\big)^{-1},
\end{split}\end{equation*}
we get, for $z \in \Gamma_{\theta,\rho}^{\overline{\tau}}$
\begin{equation*}\begin{split}
&\quad\big\|
zK(z)H(z)-
\psi_{p,\overline{\tau}}(e^{-z\overline{\tau}}) K(\psi_{p,\overline{\tau}}(e^{-z\overline{\tau}}))H_{\overline{\tau}}(z)
\big\|
\\&
\lesssim
\big\|(z^\alpha-\Delta_h)^{-1}\big\|
\big|\psi_{p,\overline{\tau}}(e^{-z\overline{\tau}})^\alpha-z^\alpha\big|
\big\|
(\psi_{p,\overline{\tau}}(e^{-z\overline{\tau}})^\alpha-\Delta_h)^{-1}
\big\|
|H(z)|
\\&\quad
+
\big\|
(\psi_{p,\overline{\tau}}(e^{-z\overline{\tau}})^\alpha-\Delta_h)^{-1}
\big\|
|H(z)-H_{\overline{\tau}}(z)|,
\end{split}\end{equation*}
which combining (\ref{App.3.1}), (\ref{App.5.1.1}) with the estimates in Lemma \ref{lem.0}, leads to the result (\ref{App.14.0.1}) required.
\end{proof}
\end{lemma}
\begin{lemma}\label{lem.0.2}
Assume functions $\mathcal{H}_{\overline{\tau}}(z)$ and $\mathcal{H}(z)$ satisfy $\|\mathcal{H}_{\overline{\tau}}(z)\|\lesssim \overline{\tau}^{s_1}|z|^{s_2}$ for $z \in \Gamma_{\theta,\rho}^{\overline{\tau}}$
and
$\|\mathcal{H}(z)\|\lesssim |z|^{s_3}$ for $z \in \Gamma_{\theta,\rho}$.
$s_1, s_2$ and $s_3$ denote some contants.
There hold
\begin{equation}\begin{split}\label{App.14.0.2}
\bigg\|
\int_{\Gamma_{\theta,\rho}^{\overline{\tau}}}e^{z\overline{t}_n}\mathcal{H}_{\overline{\tau}}(z)
\mathrm{d}z
\bigg\|
\lesssim
\overline{\tau}^{s_1}\overline{t}_n^{-s_2-1},
\quad
\bigg\|
\int_{\Gamma_{\theta,\rho} \setminus \Gamma_{\theta,\rho}^{\overline{\tau}}}e^{z\overline{t}_n}
\mathcal{H}(z)
\mathrm{d}z
\bigg\|
\lesssim
\overline{\tau}^p
\overline{t}_n^{-s_3-p-1}.
\end{split}\end{equation}
\begin{proof}
Since the contour $\Gamma_{\theta,\rho}^{\overline{\tau}}$ consists of three parts
\begin{equation*}\begin{split}
\Gamma_{\theta,\rho}^{\overline{\tau}}
=\{z=re^{\pm{\rm i}\theta}: \rho \leq r \leq \pi/(\overline{\tau}\sin\theta)\}
\cup
\{z=\rho e^{{\rm i}\sigma}: -\theta \leq \sigma \leq \theta\},
\end{split}\end{equation*}
by symmetry we obtain
\begin{equation*}\begin{split}
\bigg\|
\int_{\Gamma_{\theta,\rho}^{\overline{\tau}}}e^{z\overline{t}_n}\mathcal{H}_{\overline{\tau}}(z)
\mathrm{d}z
\bigg\|
&\lesssim
\int_\rho^{\frac{\pi}{\overline{\tau}\sin\theta}}
e^{r\overline{t}_n\cos\theta}
\|\mathcal{H}_{\overline{\tau}}(re^{{\rm i}\theta})\|
\mathrm{d}r
+\int_0^\theta
e^{\rho \overline{t}_n \cos\sigma}
\|\mathcal{H}_{\overline{\tau}}(\rho e^{{\rm i}\sigma})\|
\rho\mathrm{d}\sigma
\\
&\lesssim
\overline{\tau}^{s_1}
\bigg(
\int_\rho^{\frac{\pi}{\overline{\tau}\sin\theta}}
e^{r\overline{t}_n\cos\theta}
r^{s_2}
\mathrm{d}r
+
\rho^{s_2+1}
\int_0^\theta
e^{\rho \overline{t}_n \cos\sigma}
\mathrm{d}\sigma
\bigg).
\end{split}\end{equation*}
By choosing $\rho =\overline{t}_n^{-1}$ and setting $r\overline{t}_n=s$, one gets
\begin{equation*}\begin{split}
\int_\rho^{\frac{\pi}{\overline{\tau}\sin\theta}}
e^{r\overline{t}_n\cos\theta}
r^{s_2}
\mathrm{d}r
&\leq
\overline{t}_n^{-s_2-1}
\int_1^\infty
e^{s\cos\theta}
s^{s_2}
\mathrm{d}s
\lesssim
\overline{t}_n^{-s_2-1},
\\
\rho^{s_2+1}
\int_0^\theta
e^{\rho \overline{t}_n \cos\sigma}
\mathrm{d}\sigma
&\lesssim
\overline{t}_n^{-s_2-1},
\end{split}\end{equation*}
which yields the first result.
\par
By the symmetry of the contour, we have
\begin{equation*}\begin{split}
\bigg\|
\int_{\Gamma_{\theta,\rho} \setminus \Gamma_{\theta,\rho}^{\overline{\tau}}}e^{z\overline{t}_n}
\mathcal{H}(z)
\mathrm{d}z
\bigg\|
\lesssim
\int_{\frac{\pi}{\overline{\tau}\sin\theta}}^\infty
e^{r\overline{t}_n\cos\theta}
r^{s_3}
\mathrm{d}r
\lesssim
\overline{\tau}^p
\int_{\frac{\pi}{\overline{\tau}\sin\theta}}^\infty
e^{r\overline{t}_n\cos\theta}
r^{s_3+p}
\mathrm{d}r,
\end{split}\end{equation*}
where in the last estimate we have used the fact $r\geq \frac{\pi}{\overline{\tau}\sin\theta}$.
By setting $r\overline{t}_n=s$, one can readily obtain the second result of (\ref{App.14.0.2}).
\end{proof}
\end{lemma}

\begin{lemma}\label{lem.1}
Assume the coefficients $b_n^{(p)}$ and $d_{j,n}^{(p,\beta)}$ validate the estimates in (\ref{App.5.2}).
For $z \in \Gamma_{\theta,\rho}^{\overline{\tau}}$, there hold
\begin{itemize}
\item[$\mathrm{(i)}$]
$
\big\|K(z)-K(\psi_{p,\overline{\tau}}(e^{-z\overline{\tau}})) \mu(e^{-z\overline{\tau}})\big\|
\lesssim \overline{\tau}^p|z|^{p-1-\alpha},
$
\item[$\mathrm{(ii)}$]
$\displaystyle
\bigg\|\frac{1}{j!}\Gamma(j+1+\beta)K(z)z^{-j-\beta}
-\psi_{p,\overline{\tau}}(e^{-z\overline{\tau}})K(\psi_{p,\overline{\tau}}(e^{-z\overline{\tau}}))\overline{\tau}^{j+1+\beta}
\beta_j(e^{-z\overline{\tau}})\bigg\|
\lesssim
\overline{\tau}^p |z|^{p-j-1-\alpha-\beta}
$, with $j=0,1,2,\cdots,p-2$,
\item[$\mathrm{(iii)}$]
$\displaystyle
\bigg\|\Gamma(p+\beta)K(z)z^{-p+1-\beta}
-\psi_{p,\overline{\tau}}(e^{-z\overline{\tau}})K(\psi_{p,\overline{\tau}}(e^{-z\overline{\tau}}))\overline{\tau}^{p+\beta}
\mathrm{Li}_{-(p-1+\beta)}(e^{-z\overline{\tau}})\bigg\|
\lesssim
\overline{\tau}^p |z|^{-\alpha-\beta}.
$
\end{itemize}
\begin{proof}
The argument for (i) can be found in Lemma C.1 in \cite{jin2017correction} and we next give the details of the estimate of (ii) and (iii).
\par
For (ii), let $H(z)=\frac{1}{j!}\Gamma(j+1+\beta)z^{-j-\beta-1}$ and $H_{\overline{\tau}}(z)=\overline{\tau}^{j+1+\beta}\beta_j(e^{-z\overline{\tau}})$.
Using the second estimate in (\ref{App.5.2}) and Remark \ref{rem.1}, one gets
\begin{equation*}\begin{split}
|H(z)-H_{\overline{\tau}}(z)|
&=
\overline{\tau}^{j+1+\beta}
\bigg|
\frac{1}{j!}\Gamma(j+1+\beta)
(z\overline{\tau})^{-j-\beta-1}
-\beta_j(e^{-z\overline{\tau}})
\bigg|
\\
&\lesssim
\overline{\tau}^{j+1+\beta}
\bigg|
\beta_j(e^{-z\overline{\tau}})-
\frac{1}{j!}\Gamma(j+1+\beta)
\psi_p(e^{-z\overline{\tau}})^{-j-1-\beta}
\bigg|
\\&\quad+
\bigg|
\psi_{p,\overline{\tau}}(e^{-z\overline{\tau}})^{-j-1-\beta}-z^{-j-\beta-1}
\bigg|
\\&\lesssim
\overline{\tau}^{j+1+\beta}|1-e^{-z\overline{\tau}}|^{p-j-1-\beta}
+\overline{\tau}^p |z|^{p-j-1-\beta}
\\&\lesssim
\overline{\tau}^p |z|^{p-j-1-\beta}.
\end{split}\end{equation*}
The result of (ii) then follows from Lemma \ref{lem.0.1} and the above estimate.
\par
For (iii), set $H(z)=\Gamma(p+\beta)z^{-p-\beta}$ and $H_{\overline{\tau}}(z)=\overline{\tau}^{p+\beta}{\rm Li}_{-(p-1+\beta)}(e^{-z\overline{\tau}})$.
In accordance to Lemma \ref{lem.0} and the estimate in (\ref{App.5.4}) with $j=p-1$ and $\zeta=e^{-z\overline{\tau}}$, we obtain
\begin{equation*}\begin{split}
|H(z)-H_{\overline{\tau}}(z)|
&\lesssim
\overline{\tau}^{p+\beta}
\bigg|
{\rm Li}_{-(p-1+\beta)}(e^{-z\overline{\tau}})
-\frac{\Gamma(p+\beta)}{\psi_p(e^{-z\overline{\tau}})^{p+\beta}}
\bigg|
\\&\lesssim
\overline{\tau}^{p+\beta}
|1-e^{-z\overline{\tau}}|^{-\beta}
\lesssim
\overline{\tau}^p |z|^{-\beta},
\end{split}\end{equation*}
which leads to (iii) by resorting to Lemma \ref{lem.0.1}.
The proof for this lemma is completed.
\end{proof}
\end{lemma}
\begin{lemma}\label{lem.2}
Suppose $\phi(t)=(\log\frac{t}{a})^\beta [(\log \frac{t}{a})^{p-1} \ast \phi_1(t)]$ where $\beta \in (0,1)$ and $\phi_1(t)$ is sufficiently smooth with respect to $t$.
Then, there hold $\delta_t^j\phi(a)=0$ for $j=0,1,2,\cdots,p-1$ and
\begin{equation*}\begin{split}
|\delta_t^p \phi(t)|
\lesssim
\bigg(\log\frac{t}{a}\bigg)^{\beta-1}
\int_a^t |\phi_1(s)|
\frac{\mathrm{d}s}{s}
+
\bigg(\log\frac{t}{a}\bigg)^{\beta}
|\phi_1(t)|.
\end{split}\end{equation*}
\begin{proof}
Set $t=ae^{\overline{t}}$ and define $\Phi(\overline{t}):=\phi(t)$ and $\Phi_1(\overline{t}):=\phi_1(t)$.
In accordance to (\ref{Pre.1.2}), one gets $\Phi(\overline{t})=\overline{t}^\beta[\overline{t}^{p-1}\overline{\ast}\Phi_1(\overline{t})]$.
By using the general Leibniz rule, we have
\begin{equation}\begin{split}\label{App.14.0.3}
\Phi^{(j)}(\overline{t})
&=
\sum_{k=0}^j
\binom{j}{k}
(\overline{t}^\beta)^{(j-k)}
[\overline{t}^{p-1} 
\overline{\ast} \Phi_1(\overline{t})]^{(k)}
\\&=
\sum_{k=0}^j
\binom{j}{k}
\frac{(p-1)!\Gamma(\beta+1)}{(p-k-1)!\Gamma(\beta-j+k+1)}
\overline{t}^{\beta-j+k}
\big[
\overline{t}^{p-k-1}\overline{\ast}\Phi_1(\overline{t})
\big],
\end{split}\end{equation}
which leads to the estimate
\begin{equation*}\begin{split}
|\Phi^{(j)}(\overline{t})|
\lesssim
\overline{t}^{p+\beta-j-1}
\int_0^{\overline{t}}|\Phi_1(s)|
\mathrm{d}s,\quad
j=0,1,\cdots,p-1.
\end{split}\end{equation*}
Combining the fact $\delta_t^j\phi(a)=\Phi^{(j)}(0)$ with the above estimate, we get
$\delta_t^j\phi(a)=0$ for $j=0,1,2,\cdots,p-1$.
Let $\ell_k=\frac{(p-1)!\Gamma(\beta+1)}{(p-k-1)!\Gamma(\beta-j+k+1)}$ and rewrite (\ref{App.14.0.3}) with $j=p-1$ to yield
\begin{equation*}\begin{split}
\Phi^{(p-1)}(\overline{t})
=
\sum_{k=0}^{p-2}\ell_k
\overline{t}^{\beta-p+1+k}
[\overline{t}^{p-k-1} \overline{\ast} \Phi_1(\overline{t})]
+(p-1)!\overline{t}^\beta
[1\overline{\ast}\Phi_1(\overline{t})].
\end{split}\end{equation*}
Taking the first derivative of the above function, one gets
\begin{equation*}\begin{split}
\Phi^{(p)}(\overline{t})
&=
\sum_{k=0}^{p-2}
\bigg\{
(\beta-p+1+k)\ell_k
\overline{t}^{\beta-p+k}
[\overline{t}^{p-k-1} \overline{\ast} \Phi_1(\overline{t})]
+(p-k-1)\ell_k
\overline{t}^{\beta-p+1+k}
[\overline{t}^{p-k-2} \overline{\ast} \Phi_1(\overline{t})]
\bigg\}
\\&\quad+
(p-1)!
\bigg\{
\beta \overline{t}^{\beta-1}
[1\overline{\ast}\Phi_1(\overline{t})]
+\overline{t}^\beta \Phi_1(\overline{t})
\bigg\},
\end{split}\end{equation*}
which indicates the following estimate
\begin{equation}\begin{split}\label{App.14.0.4}
|\Phi^{(p)}(\overline{t})|
\lesssim
\overline{t}^{\beta-1}
\int_0^{\overline{t}}
|\Phi_1(s)|
\mathrm{d}s
+\overline{t}^\beta
|\Phi_1(\overline{t})|.
\end{split}\end{equation}
Replacing $\overline{t}$ with $\log\frac{t}{a}$ and using the fact $\delta_t^p\phi(t)=\Phi^{(p)}(\overline{t})$, we complete the proof of the lemma.
\end{proof}
\end{lemma}

\par
Drawing on the results of Steps 1 to 3 above, we present the sharp error estimates for the corrected high-order scheme (\ref{App.5.1}) in the next several theorems.
\begin{theorem}\label{thm.3}
Assume $v(\boldsymbol x) \equiv 0, \beta \in [0,1)$ and $f(\boldsymbol x,t)\equiv \big(\log\frac{t}{a}\big)^{p-1+\beta}f_1(\boldsymbol x)$.
Let $w_h(t)$ be the solutions of (\ref{App.5}) and $W_h^n$ be the solution of (\ref{App.5.1}) without correction terms.
There holds
\begin{equation*}\begin{split}
\|w_h(t_n)-W_h^n\|\lesssim
\overline{\tau}^p
\bigg(\log\frac{t_n}{a}\bigg)^{\alpha+\beta-1}\|f_1(\boldsymbol x)\|.
\end{split}\end{equation*}
\begin{proof}
With the assumption $v(\boldsymbol x) \equiv 0$ and $f(\boldsymbol x,t)\equiv \big(\log\frac{t}{a}\big)^{p-1+\beta}f_1(\boldsymbol x)$, one gets from (\ref{App.7.1}) that
\begin{equation}\begin{split} \label{App.14.1}
w_h(t)
&=\frac{1}{2\pi\rm{i}}\int_{\Gamma_{\theta,\rho}}e^{z\overline{t}}zK(z)\widehat{R}_p(z)\mathrm{d}z,
\quad
\overline{t}=\log\frac{t}{a},
\end{split}\end{equation}
where $\widehat{R}_p(z)$ is defined by
\begin{equation*}\begin{split}
\widehat{R}_p(z)=f_1(\boldsymbol x)\mathcal{L}\bigg\{\bigg(\log\frac{t}{a}\bigg)^{p-1+\beta}\bigg\}(z)
=f_1(\boldsymbol x)\frac{\Gamma(p+\beta)}{z^{p+\beta}}.
\end{split}\end{equation*}
In accordance with (\ref{App.14}), we have
\begin{equation}\begin{split}\label{App.15}
W_h^n=\frac{1}{2\pi\rm{i}}\int_{\Gamma_{\theta,\rho}^{\overline{\tau}}}
e^{z\overline{t}_n}\psi_{p,\overline{\tau}}(e^{-z\overline{\tau}}) K(\psi_{p,\overline{\tau}}(e^{-z\overline{\tau}}))
\overline{\tau}\widetilde{R}_p(e^{-z\overline{\tau}})
\mathrm{d}z,
\end{split}\end{equation}
where
\begin{equation*}\begin{split}
\widetilde{R}_p(\zeta)=
\sum_{n=0}^\infty \zeta^n f(\boldsymbol x,t_n)
=f_1(\boldsymbol x)\sum_{n=0}^\infty \zeta^n\bigg(\log\frac{t_n}{a}\bigg)^{p-1+\beta}
=f_1(\boldsymbol x)\overline{\tau}^{p-1+\beta}\mathrm{Li}_{-(p-1+\beta)}(\zeta).
\end{split}\end{equation*}
Then, by (\ref{App.14.1}) and (\ref{App.15}), one immediately gets
\begin{equation*}\begin{split}
w_h(t_n)-W_h^n
=J+J',
\end{split}\end{equation*}
where
\begin{equation*}\begin{split}
J:&
=\frac{1}{2\pi\rm{i}}\int_{\Gamma_{\theta,\rho}^{\overline{\tau}}}e^{z\overline{t}_n}
\big(zK(z)\widehat{R}_p(z)
-\psi_{p,\overline{\tau}}(e^{-z\overline{\tau}}) K(\psi_{p,\overline{\tau}}(e^{-z\overline{\tau}}))
\overline{\tau}\widetilde{R}_p(e^{-z\overline{\tau}})\big)\mathrm{d}z,
\\
J':&
=\frac{1}{2\pi\rm{i}}\int_{\Gamma_{\theta,\rho} \setminus \Gamma_{\theta,\rho}^{\overline{\tau}}}e^{z\overline{t}_n}zK(z)\widehat{R}_p(z)\mathrm{d}z.
\end{split}\end{equation*}
Using the result (iii) in Lemma \ref{lem.1} and Lemma \ref{lem.0.2}, we have
\begin{equation*}\begin{split}
\|J\|+\|J'\|
&\lesssim
\overline{\tau}^p\overline{t}_n^{\alpha+\beta-1}\|f_1(\boldsymbol x)\|,
\end{split}\end{equation*}
which, by replacing $\overline{t}_n$ with $\log\frac{t_n}{a}$, completes the proof of the theorem.
\end{proof}
\end{theorem}
\begin{remark}\label{rem.2}
It is notable that the result of the above theorem can be extended to the continuous case
\begin{equation}\begin{split}\label{App.16}
\|w_h(t)-W_h(t)\|\lesssim
\overline{\tau}^p
\bigg(\log\frac{t}{a}\bigg)^{\alpha+\beta-1}\|f_1(\boldsymbol x)\|,\quad
t\in(a,T],
\end{split}\end{equation}
where $W_h(t)$ should be understood in the following way that by introducing the discrete solution operator
\begin{equation*}\begin{split}
\mathcal{E}_{\overline{\tau}}^j=
\frac{\overline{\tau}}{2\pi{\rm i}}\int_{\Gamma_{\theta,\rho}^{\overline{\tau}}}e^{z\overline{t}_j}
\psi_{p,\overline{\tau}}(e^{-z\overline{\tau}})
K\big(\psi_{p,\overline{\tau}}(e^{-z\overline{\tau}})\big)
\mathrm{d}z,
\end{split}\end{equation*}
and defining $\mathcal{E}_{\overline{\tau}}(\overline{t}):=\sum_{j=0}^\infty \mathcal{E}_{\overline{\tau}}^j \Delta_{\overline{t}_j}(\overline{t})$, where $\Delta_{\overline{t}_j}(\overline{t})$ is the Dirac delta function at $\overline{t}_j$,  for the source term $F(\overline{t}):=f(\boldsymbol x,t)$, we can obtain $W_h(t)$ such that
\begin{equation*}\begin{split}
W_h(t)=[\mathcal{E}_{\overline{\tau}}(\overline{t}) \overline{\ast} F(\overline{t})]
\big|_{\overline{t}=\log\frac{t}{a}}.
\end{split}\end{equation*}
Similarly, by introducing the continuous solution operator $\mathcal{E}(\overline{t})=\frac{1}{2\pi \mathrm{i}}\int_{\Gamma_{\theta,\rho}}e^{z\overline{t}}zK(z)\mathrm{d}z$, the solution $w_h(t)$ can be derived by
\begin{equation*}\begin{split}
w_h(t)=[\mathcal{E}(\overline{t}) \overline{\ast} F(\overline{t})]
\big|_{\overline{t}=\log\frac{t}{a}}.
\end{split}\end{equation*}
Therefore, Theorem \ref{thm.3} actually indicates that
\begin{equation*}\begin{split}
\bigg\|
\big(
\mathcal{E}_{\overline{\tau}}(\overline{t})
-
\mathcal{E}(\overline{t})
\big)
\overline{\ast}
F(\overline{t})
\bigg\|_{\overline{t}=\log\frac{t_n}{a}}
\lesssim
\overline{\tau}^p
\bigg(\log\frac{t_n}{a}\bigg)^{\alpha+\beta-1}\|f_1(\boldsymbol x)\|,
\end{split}\end{equation*}
which can be generalized to the continuous case (\ref{App.16}) by similar technique presented in Lemma 3.12 in \cite{jin2018analysis}.
\end{remark}

\begin{theorem}\label{thm.4}
Assume $v(\boldsymbol x) \equiv 0, \beta\in[0,1)$ and $f(\boldsymbol x,t)\equiv \big(\log\frac{t}{a}\big)^{\beta}\big[\big(\log\frac{t}{a}\big)^{p-1}\ast f_2(\boldsymbol x,t)\big]$ where $f_2$ is smooth with respect to $t$.
Let $w_h(t)$ be the solutions of (\ref{App.5}) and $W_h^n$ be the solution of (\ref{App.5.1}) without correction terms.
Then, there holds
\begin{equation}\begin{split}\label{App.16.1}
\|w_h(t_n)-W_h^n\|
&\lesssim
\overline{\tau}^p 
\bigg[
\kappa_\beta
\bigg(\log\frac{t_n}{a}\bigg)^{\alpha+\beta-1}
\int_a^{t_n}
\|f_2(\boldsymbol x,s)\|\frac{\mathrm{d}s}{s}
\\&\quad\quad\quad
+
\int_a^{t_n}
\bigg(\log\frac{t_n}{s}\bigg)^{\alpha-1}
\bigg(\log\frac{s}{a}\bigg)^\beta
\|f_2(\boldsymbol x,s)\|
\frac{\mathrm{d}s}{s}
\bigg],
\end{split}\end{equation}
where $\kappa_\beta$ is an indicator satisfying $\kappa_\beta=0$ for $\beta=0$ and $\kappa_\beta=1$ for $\beta \in (0,1)$. 
\begin{proof}
For given $t$, define $\overline{t}$ such that $t=ae^{\overline{t}}$ and take $F(\boldsymbol x,\overline{t}):=f(\boldsymbol x,t)$ and $F_2(\boldsymbol x,\overline{t}):=f_2(\boldsymbol x,t)$.
It holds that $F(\boldsymbol x,\overline{t})=\overline{t}^\beta [\overline{t}^{p-1} \overline{\ast} F_2(\boldsymbol x,\overline{t})]$.
\par
\textbf{Case I: } $\beta\in(0,1)$.
Using the modified Taylor theorem \ref{thm.0} and the results in Lemma \ref{lem.2}, we obtain
\begin{equation*}\begin{split}
f(\boldsymbol x,t)=\frac{1}{(p-1)!}
(\log\frac{t}{a})^{p-1}
\ast
\delta_t^p f(\boldsymbol x,t)
=\frac{1}{(p-1)!}
\overline{t}^{p-1}
\overline{\ast}
F^{(p)}(\boldsymbol x,\overline{t}).
\end{split}\end{equation*}
By resorting to the continuous and discrete solution operators introduced in Remark \ref{rem.2}, one gets
\begin{equation*}\begin{split}
w_h(t_n)-W_h^n
=
\frac{1}{(p-1)!}
\big[
\big(\mathcal{E}(\overline{t})-\mathcal{E}_{\overline{\tau}}(\overline{t})\big)
\overline{\ast}
\overline{t}^{p-1}
\overline{\ast}
F^{(p)}(\boldsymbol x,\overline{t})
\big]|\big|_{\overline{t}=\log\frac{t_n}{a}},
\end{split}\end{equation*}
implying the estimate
\begin{equation}\begin{split}\label{App.17}
\|w_h(t_n)-W_h^n\|
\lesssim
\overline{\tau}^p
\int_0^{\overline{t}_n}
(\overline{t}_n-s)^{\alpha-1}
\|F^{(p)}(\boldsymbol x,s)\|
\mathrm{d}s,\quad
\overline{t}_n=\log\frac{t_n}{a}.
\end{split}\end{equation}
Recalling that (\ref{App.14.0.4}) indicates the estimate
\begin{equation*}\begin{split}
\|F^{(p)}(\boldsymbol x,s)\|
\lesssim
s^{\beta-1}
\int_0^s
\|F_2(\boldsymbol x,\ell)\|
\mathrm{d}\ell
+s^\beta
\|F_2(\boldsymbol x,s)\|,
\end{split}\end{equation*}
which, in combination with (\ref{App.17}), leads to
\begin{equation*}\begin{split}
\|w_h(t_n)-W_h^n\|
&\lesssim
\overline{\tau}^p
\bigg[
\int_0^{\overline{t}_n}
(\overline{t}_n-s)^{\alpha-1}
s^{\beta-1}
\int_0^s
\|F_2(\boldsymbol x,\ell)\|
\mathrm{d}\ell
\mathrm{d}s
+
\int_0^{\overline{t}_n}
(\overline{t}_n-s)^{\alpha-1}
s^{\beta}
\|F_2(\boldsymbol x,s)\|
\mathrm{d}s
\bigg]
\\&\lesssim
\overline{\tau}^p
\bigg[
\overline{t}_n^{\alpha+\beta-1}
\int_0^{\overline{t}_n}
\|F_2(\boldsymbol x,\ell)\|
\mathrm{d}\ell
+
\int_0^{\overline{t}_n}
(\overline{t}_n-s)^{\alpha-1}
s^{\beta}
\|F_2(\boldsymbol x,s)\|
\mathrm{d}s
\bigg].
\end{split}\end{equation*}
By replacing $\overline{t}_n$ with $\log\frac{t_n}{a}$, we get the result (\ref{App.16.1}) for $\beta \in (0,1)$.
\par
\textbf{Case II: } $\beta=0$.
In this case, since $F(\boldsymbol x,\overline{t})=\overline{t}^{p-1} \overline{\ast} F_2(\boldsymbol x,\overline{t})$ is of nonsingular, the error estimate can be carried out by similar argument as that in Lemma 3.12 in \cite{jin2018analysis}, which will be omitted here for space reasons.
We note that in this case the error estimate is of the form
\begin{equation*}\begin{split}
\|w_h(t_n)-W_h^n\|
&\lesssim
\overline{\tau}^p
\int_0^{\overline{t}_n}
(\overline{t}_n-s)^{\alpha-1}
\|F_2(\boldsymbol x,s)\|
\mathrm{d}s
\\&=
\overline{\tau}^p 
\int_a^{t_n}
\bigg(\log\frac{t_n}{s}\bigg)^{\alpha-1}
\|f_2(\boldsymbol x,s)\|
\frac{\mathrm{d}s}{s},
\end{split}\end{equation*}
which confirms the result (\ref{App.16.1}) and the proof of the theorem is completed.
\end{proof}
\end{theorem}
\begin{theorem}\label{thm.5}
For the source term $f(\boldsymbol x,t)$ given in (\ref{App.1.1}),
assume $w_h(t)$ is the solution of (\ref{App.5}) and $W_h^n$ the solution of (\ref{App.5.1}).
Let $U_h^n:=W_h^n+v_h$ be the approximation to $u_h(t_n):=w_h(t_n)+v_h$ with $v_h=R_h v$.
There holds the following error estimate
\begin{equation}\begin{split}\label{App.17.1}
\|u_h(t_n)-U_h^n\|
&\lesssim
\overline{\tau}^p
\bigg[
 \bigg(
\log\frac{t_n}{a}
\bigg)^{\alpha-p}
(\|\Delta v\|+\|f_h^0\|)
+\sum_{j=0}^{p-1}
\bigg(
\log\frac{t_n}{a}
\bigg)^{-p+j+\alpha+\beta}
\big\|\delta_t^j g_h(a)\big\|
\\&\quad\quad+
\kappa_\beta
\bigg(\log\frac{t_n}{a}\bigg)^{\alpha+\beta-1}
\int_a^{t_n}
\big\|
\delta_s^p g_h
\big\|
\frac{\mathrm{d}s}{s}
+
\int_a^{t_n}
\bigg(\log\frac{t_n}{s}\bigg)^{\alpha-1}
\bigg(\log\frac{s}{a}\bigg)^\beta
\big\|
\delta_s^p g_h
\big\|
\frac{\mathrm{d}s}{s}
\bigg],
\end{split}\end{equation}
where $\kappa_\beta$ is an indicator satisfying $\kappa_\beta=0$ for $\beta=0$ and $\kappa_\beta=1$ for $\beta \in (0,1)$. 
\begin{proof}
In accordance to the solution representation (\ref{App.7.1}) and (\ref{App.14}),
we split the error into several parts:
\begin{equation*}\begin{split}
u_h(t_n)-U_h^n=w_h(t_n)-W_h^n=J_1+J_1'+J_2+J_2'+J_3+J_3',
\end{split}\end{equation*}
where
\begin{equation*}\begin{split}
J_1:&=\frac{1}{2\pi\rm{i}}\int_{\Gamma_{\theta,\rho}^{\overline{\tau}}}e^{z\log\frac{t_n}{a}}K(z)(\Delta_h v_h+f_h^0)
-e^{z\overline{t}_n}K(\psi_{p,\overline{\tau}}(e^{-z\overline{\tau}})) \mu(e^{-z\overline{\tau}})(\Delta_h v_h+f_h^0)
\mathrm{d}z
\\
&=
\frac{1}{2\pi\rm{i}}\int_{\Gamma_{\theta,\rho}^{\overline{\tau}}}e^{z \overline{t}_n}
\big[K(z)-K(\psi_{p,\overline{\tau}}(e^{-z\overline{\tau}})) \mu(e^{-z\overline{\tau}})\big](\Delta_h v_h+f_h^0)
\mathrm{d}z,
\\
J_1':&=\frac{1}{2\pi\rm{i}}\int_{\Gamma_{\theta,\rho}\setminus \Gamma_{\theta,\rho}^{\overline{\tau}}}e^{z\log\frac{t_n}{a}}K(z)(\Delta_h v_h+f_h^0)\mathrm{d}z
=\frac{1}{2\pi\rm{i}}\int_{\Gamma_{\theta,\rho}\setminus \Gamma_{\theta,\rho}^{\overline{\tau}}}e^{z\overline{t}_n}K(z)(\Delta_h v_h+f_h^0)\mathrm{d}z,
\end{split}\end{equation*}
\begin{equation*}\begin{split}
J_2:&=\frac{1}{2\pi\rm{i}}\int_{\Gamma_{\theta,\rho}^{\overline{\tau}}}e^{z\log\frac{t_n}{a}}zK(z)\sum_{j=0}^{p-2}\frac{\Gamma(j+1+\beta)}{j! z^{j+1+\beta}}\delta_t^j g_h(a)\mathrm{d}z
\\
&\quad-\frac{1}{2\pi\rm{i}}\int_{\Gamma_{\theta,\rho}^{\overline{\tau}}}
e^{z\overline{t}_n}\psi_{p,\overline{\tau}}(e^{-z\overline{\tau}}) K(\psi_{p,\overline{\tau}}(e^{-z\overline{\tau}}))
\sum_{j=0}^{p-2}\overline{\tau}^{j+1+\beta} \delta_t^j g_h(a)
\beta_j(e^{-z\overline{\tau}})\mathrm{d}z
\\
&=\frac{1}{2\pi\rm{i}}\int_{\Gamma_{\theta,\rho}^{\overline{\tau}}}e^{z\overline{t}_n}
\bigg[zK(z)\sum_{j=0}^{p-2}\frac{\Gamma(j+1+\beta)}{j! z^{j+1+\beta}} \delta_t^j g_h(a)
\\&\quad\quad\quad\quad\quad\quad\quad\quad
-\psi_{p,\overline{\tau}}(e^{-z\overline{\tau}}) K(\psi_{p,\overline{\tau}}(e^{-z\overline{\tau}}))
\sum_{j=0}^{p-2}\overline{\tau}^{j+1+\beta} \delta_t^j g_h(a)
\beta_j(e^{-z\overline{\tau}})\bigg]\mathrm{d}z
\\&=
\sum_{j=0}^{p-2}\frac{1}{2\pi\rm{i}}\int_{\Gamma_{\theta,\rho}^{\overline{\tau}}}e^{z\overline{t}_n}\bigg[
\frac{1}{j!}\Gamma(j+1+\beta)K(z)z^{-j-\beta}
\\&\quad\quad\quad\quad\quad\quad\quad\quad\quad
-\psi_{p,\overline{\tau}}(e^{-z\overline{\tau}})K(\psi_{p,\overline{\tau}}(e^{-z\overline{\tau}}))\overline{\tau}^{j+1+\beta}
\beta_j(e^{-z\overline{\tau}})
\bigg]\delta_t^j g_h(a)\mathrm{d}z
\\&=: \sum_{j=0}^{p-2}J_{2,j},
\\
J_2':&=\frac{1}{2\pi\rm{i}}\int_{\Gamma_{\theta,\rho}\setminus \Gamma_{\theta,\rho}^{\overline{\tau}}}e^{z\log\frac{t_n}{a}}zK(z)\sum_{j=0}^{p-2}\frac{\Gamma(j+1+\beta)}{j! z^{j+1+\beta}}\delta_t^j g_h(a)\mathrm{d}z
\\
&=\frac{1}{2\pi\rm{i}}\int_{\Gamma_{\theta,\rho}\setminus \Gamma_{\theta,\rho}^{\overline{\tau}}}e^{z\overline{t}_n}zK(z)\sum_{j=0}^{p-2}\frac{\Gamma(j+1+\beta)}{j! z^{j+1+\beta}} \delta_t^j g_h(a)\mathrm{d}z
\\
&=
\sum_{j=0}^{p-2}
\frac{1}{2\pi\rm{i}}\int_{\Gamma_{\theta,\rho}\setminus \Gamma_{\theta,\rho}^{\overline{\tau}}}e^{z\overline{t}_n}
\frac{1}{j!}
\Gamma(j+1+\beta)
K(z)z^{-j-\beta} 
\delta_t^j g_h(a)
\mathrm{d}z
\\
&=:
\sum_{j=0}^{p-2}J'_{2,j},
\end{split}\end{equation*}
\begin{equation*}\begin{split}
J_3:&=\frac{1}{2\pi\rm{i}}\int_{\Gamma_{\theta,\rho}^{\overline{\tau}}}e^{z\log\frac{t_n}{a}}zK(z)\widehat{R}_p(z)\mathrm{d}z
-\frac{1}{2\pi\rm{i}}\int_{\Gamma_{\theta,\rho}^{\overline{\tau}}}
e^{z\overline{t}_n}\psi_{p,\overline{\tau}}(e^{-z\overline{\tau}}) K(\psi_{p,\overline{\tau}}(e^{-z\overline{\tau}}))
\overline{\tau}\widetilde{R}_p(e^{-z\overline{\tau}})\mathrm{d}z
\\
&=\frac{1}{2\pi\rm{i}}\int_{\Gamma_{\theta,\rho}^{\overline{\tau}}}e^{z\overline{t}_n}
\big(zK(z)\widehat{R}_p(z)
-\psi_{p,\overline{\tau}}(e^{-z\overline{\tau}}) K(\psi_{p,\overline{\tau}}(e^{-z\overline{\tau}}))
\overline{\tau}\widetilde{R}_p(e^{-z\overline{\tau}})\big)\mathrm{d}z,
\\
J_3':&=\frac{1}{2\pi\rm{i}}\int_{\Gamma_{\theta,\rho} \setminus \Gamma_{\theta,\rho}^{\overline{\tau}}}e^{z\log\frac{t_n}{a}}zK(z)\widehat{R}_p(z)\mathrm{d}z
=\frac{1}{2\pi\rm{i}}\int_{\Gamma_{\theta,\rho} \setminus \Gamma_{\theta,\rho}^{\overline{\tau}}}e^{z\overline{t}_n}zK(z)\widehat{R}_p(z)\mathrm{d}z.
\end{split}\end{equation*}
\par
For the estimates of $J_1$ and $J'_1$, using (\ref{App.3.1}), (i) of Lemma \ref{lem.1} and Lemma \ref{lem.0.2}, we obtain
\begin{equation}\begin{split}\label{App.18}
\|J_1\|+\|J'_1\|
&\lesssim
\overline{\tau}^p \overline{t}_n^{\alpha-p}
\|\Delta_h v_h+f_h^0\|
\lesssim
\overline{\tau}^p \overline{t}_n^{\alpha-p}
(\|\Delta_h v_h\|+\|f_h^0\|)
\\&
\lesssim
\overline{\tau}^p \overline{t}_n^{\alpha-p}
(\|\Delta v\|+\|f_h^0\|)
=
\overline{\tau}^p 
\bigg(
\log\frac{t_n}{a}
\bigg)^{\alpha-p}
(\|\Delta v\|+\|f_h^0\|),
\end{split}\end{equation}
where in the last step we have used the fact $\Delta_h R_h=P_h \Delta$.
\par
Similarly, for each $J_{2,j}$, using (\ref{App.3.1}), (ii) of Lemma \ref{lem.1} and Lemma \ref{lem.0.2}, one gets
\begin{equation}\begin{split}\label{App.19}
\|J_{2,j}\|+\|J'_{2,j}\|
\lesssim
\overline{\tau}^p 
\overline{t}_n^{-p+j+\alpha+\beta}
\big\|\delta_t^j g_h(a)\big\|
=\overline{\tau}^p 
\bigg(
\log\frac{t_n}{a}
\bigg)^{-p+j+\alpha+\beta}
\big\|\delta_t^j g_h(a)\big\|.
\end{split}\end{equation}

For the terms $J_3$ and $J'_3$, recall that $R_p(t)$ in (\ref{App.6.1}) consists of two parts $R_p^{(1)}(t)$ and $R_p^{(2)}(t)$.
Introduce $J_3^{(i)}$ and $J_3^{(i)'}$ ($i=1,2$) such that
\begin{equation*}\begin{split}
J_3^{(i)}
&=
\frac{1}{2\pi\rm{i}}\int_{\Gamma_{\theta,\rho}^{\overline{\tau}}}e^{z\overline{t}_n}
\big(zK(z)\widehat{R}^{(i)}_p(z)
-\psi_{p,\overline{\tau}}(e^{-z\overline{\tau}}) K(\psi_{p,\overline{\tau}}(e^{-z\overline{\tau}}))
\overline{\tau}\widetilde{R}^{(i)}_p(e^{-z\overline{\tau}})\big)\mathrm{d}z
\\
J_3^{(i)'}
&=
\frac{1}{2\pi\rm{i}}\int_{\Gamma_{\theta,\rho} \setminus \Gamma_{\theta,\rho}^{\overline{\tau}}}e^{z\overline{t}_n}zK(z)\widehat{R}^{(i)}_p(z)\mathrm{d}z.
\end{split}\end{equation*}
Then, Theorem \ref{thm.3} indicates that
\begin{equation}\begin{split}\label{App.20}
\|J_3^{(1)}\|+\|J_3^{(1)'}\|
\lesssim
\overline{\tau}^p
\bigg(\log\frac{t_n}{a}\bigg)^{\alpha+\beta-1}
\|\delta_t^{p-1}g_h(a)\|,
\end{split}\end{equation}
and, Theorem \ref{thm.4} implies that
\begin{equation}\begin{split}\label{App.21}
\|J_3^{(2)}\|+\|J_3^{(2)'}\|
&\lesssim
\overline{\tau}^p 
\bigg[
\kappa_\beta
\bigg(\log\frac{t_n}{a}\bigg)^{\alpha+\beta-1}
\int_a^{t_n}
\big\|
\delta_s^p g_h
\big\|
\frac{\mathrm{d}s}{s}
\\&\quad\quad\quad
+
\int_a^{t_n}
\bigg(\log\frac{t_n}{s}\bigg)^{\alpha-1}
\bigg(\log\frac{s}{a}\bigg)^\beta
\big\|
\delta_s^p g_h
\big\|
\frac{\mathrm{d}s}{s}
\bigg].
\end{split}\end{equation}
The final result (\ref{App.17.1}) then follows from (\ref{App.18}), (\ref{App.19}),(\ref{App.20}) and (\ref{App.21}).
The proof is completed.
\end{proof}
\end{theorem}
\begin{remark}\label{rem.3}
We remark that for nonsmooth initial conditions such as $v(\boldsymbol x) \in L^2(\Omega)$, by setting $v_h(\boldsymbol x):=P_h v(\boldsymbol x)$, the error estimate can be presented similarly as that in Theorem \ref{thm.5} with minimal modifications (see Theorem 3.8 in \cite{jin2018analysis}), which takes the form
\begin{equation*}\begin{split}
\|u_h(t_n)-U_h^n\|
&\lesssim
\overline{\tau}^p
\bigg[
 \bigg(
\log\frac{t_n}{a}
\bigg)^{-p}
\|v\|
+ \bigg(
\log\frac{t_n}{a}
\bigg)^{\alpha-p}
\|f_h^0\|
+\sum_{j=0}^{p-1}
\bigg(
\log\frac{t_n}{a}
\bigg)^{-p+j+\alpha+\beta}
\big\|\delta_t^j g_h(a)\big\|
\\&\quad\quad+
\kappa_\beta
\bigg(\log\frac{t_n}{a}\bigg)^{\alpha+\beta-1}
\int_a^{t_n}
\big\|
\delta_s^p g_h
\big\|
\frac{\mathrm{d}s}{s}
+
\int_a^{t_n}
\bigg(\log\frac{t_n}{s}\bigg)^{\alpha-1}
\bigg(\log\frac{s}{a}\bigg)^\beta
\big\|
\delta_s^p g_h
\big\|
\frac{\mathrm{d}s}{s}
\bigg].
\end{split}\end{equation*}
It is also notable that by similar analysis one can show that the scheme (\ref{App.5.1}) is only first order accuracy (at fixed time) without any correction terms.
See Table \ref{tab2} in the next section.
\end{remark}
\section{Numerical tests}\label{sec.num}
In this section, we confirm our theoretical results for the CH fractional subdiffusion problem (\ref{App.1}) by numerical examples.
In general, high accuracy spacial discretization method is required to obtain practically the $p$th-order ($1\leq p \leq 6$) temporal accuracy of our modified scheme and therefore only the one dimensional case is considered.
We shall adopt the $P_5$ element in spacial direction.
Let $\Omega=(0,\pi)$, $t\in[1,e^2]$ and assume
$$v(x)=\sin x,\quad u(x,t)=\big[1+(\log t)^\alpha\big]\sin x,$$ such that the singular source term (with $\beta=\alpha$) is as follows
\begin{equation*}\begin{split}
f(x,t)=\big[1+\Gamma(\alpha+1)\big]\sin x
+(\log t)^\alpha 
\sin x. 
\end{split}\end{equation*}
In Table \ref{tab2}, we choose $h=1/100$ and $\alpha=0.5$.
Clearly, the fractional BDF-1 scheme needs not to be modified.
For fractional BDF-$p$ with $2 \leq p \leq 6$, one observes that the error is much smaller than that of the corresponding scheme without corrections (the standard scheme).
For cases $2 \leq p \leq 5$, by adding correction, the desired optimal convergence order can be obtained while for the case $p=6$, the numerical result is optimal than our conclusion which deserves more future exploration.
Clearly, due to the weak singularity of the solution and source term, the error is even insensitive to the fractional BDF-$p$ adopted if correction terms are omitted and only first-order accuracy is arrived at.
\begin{table}[]
\centering
\caption{Comparisons of the modified scheme (\ref{App.5.1}) with the standard scheme}\label{tab2}
{\renewcommand{\arraystretch}{1.4}
\begin{tabular*}{\hsize}{@{}@{\extracolsep{\fill}}cccclcc@{}}
\hline
\multirow{2}{*}{BDF-$p$} & \multirow{2}{*}{$\overline{\tau}$} & \multicolumn{2}{c}{Modified Scheme} &  & \multicolumn{2}{c}{Standard Scheme} \\ \cline{3-4} \cline{6-7} 
                         &                         & Error               & Order         &  & Error               & Order         \\ \hline
\multirow{3}{*}{$p=1$}   & 1/40                    & --                  & --            &  & 9.3088E-04          & --            \\
                         & 1/80                    & --                  & --            &  & 4.2987E-04          & 1.11          \\
                         & 1/160                   & --                  & --            &  & 2.0259E-04          & 1.09          \\ \hline
\multirow{3}{*}{$p=2$}   & 1/40                    & 3.0232E-05          & --            &  & 1.6418E-03          & --            \\
                         & 1/80                    & 6.9577E-06          & 2.12          &  & 7.8103E-04          & 1.07          \\
                         & 1/160                   & 1.6407E-06          & 2.08          &  & 3.7710E-04          & 1.05          \\ \hline
\multirow{3}{*}{$p=3$}   & 1/40                    & 1.7289E-06          & --            &  & 1.6254E-03          & --            \\
                         & 1/80                    & 1.9471E-07          & 3.15          &  & 7.7709E-04          & 1.06          \\
                         & 1/160                   & 2.2753E-08          & 3.10          &  & 3.7613E-04          & 1.05          \\ \hline
\multirow{3}{*}{$p=4$}   & 1/40                    & 1.4373E-07          & --            &  & 1.6262E-03          & --            \\
                         & 1/80                    & 7.8090E-09          & 4.20          &  & 7.7718E-04          & 1.07          \\
                         & 1/160                   & 4.4926E-10          & 4.12          &  & 3.7614E-04          & 1.05          \\ \hline
\multirow{3}{*}{$p=5$}   & 1/40                    & 2.8425E-08          & --            &  & 1.6261E-03          & --            \\
                         & 1/80                    & 4.1143E-10          & 6.11          &  & 7.7717E-04          & 1.07          \\
                         & 1/160                   & 1.1727E-11          & 5.13          &  & 3.7614E-04          & 1.05          \\ \hline
\multirow{3}{*}{$p=6$}   & 1/60                    & 1.1506E-06          & --            &  & 1.2639E-03          & --            \\
                         & 1/120                   & 1.0035E-10          & 13.49         &  & 6.1687E-04          & 1.03          \\
                         & 1/240                   & 1.8371E-13          & 9.09          &  & 3.0303E-04          & 1.03          \\ \hline
\end{tabular*}
}
\end{table}

\section{Conclusions}\label{sec.con}
The Hadamard fractional calculus operators are approximated by the convolution quadrature method in this work.
The local truncation error is derived with respect to the regularity of the solution.
To improve the low accuracy incurred by the weak regularity of the solution and source term of the Caputo-Hadamard fractional subdiffusion problem, novel correction technique is developed which generalizes the classical correction method designed for smooth source terms.
Rigorous and sharp error analysis is carried out in detail and is confirmed by some numerical tests.

\section{Declaration of competing interest}
The authors declare that they have no known competing financial interests or personal relationships that could have appeared to influence the work reported in this paper.
\section{Data availability}
Data will be made available on request.
\section{Acknowledgments}
This work is supported by the Research Seed Funding for High-level Talents in 2022 (No. 10000-22311201/018 to B.Y.)
and National Natural Science Foundation of China (No. 12201322 to B.Y., 12061053 to Y.L.
and 12161063 to H.L.) and Natural Science Foundation of Inner Mongolia (No. 2021BS01003 to G.Z., 2020MS01003 to Y.L., and 2021MS01018 to H.L.).

\bibliography{mybibfile}
\end{document}